\input amstex 
\input amsppt.sty 
\hsize 30pc
\vsize 47pc
\def\nmb#1#2{#2}         
\def\cit#1#2{\ifx#1!\cite{#2}\else#2\fi}
\def\totoc{}             
\def\idx{}               
\def\ign#1{}             

\redefine\o{\circ}
\define\X{\frak X}
\define\al{\alpha}
\define\be{\beta}
\define\ga{\gamma}

\define\ze{\zeta}
\define\et{\eta}
\define\th{\theta}

\define\ka{\kappa}
\define\la{\lambda}
\define\rh{\rho}
\define\si{\sigma}

\define\ph{\varphi}
\define\ch{\chi}
\define\ps{\psi}
\define\om{\omega}
\define\Ga{\Gamma}

\define\La{\Lambda}

\define\Ph{\Phi}
\define\Ps{\Psi}
\define\Om{\Omega}
\redefine\i{^{-1}}
\define\row#1#2#3{#1_{#2},\ldots,#1_{#3}}
\define\x{\times}
\define\End{\operatorname{End}}
\define\Fl{\operatorname{Fl}}

\define\sign{\operatorname{sign}}
\define\Ad{\operatorname{Ad}}
\define\ad{\operatorname{ad}}
\define\ev{\operatorname{ev}}
\redefine\L{{\Cal L}}
\define\ddt{\left.\tfrac \partial{\partial t}\right\vert_0}
\define\g{{\frak g}}
\define\h{{\frak h}}
\define\Pt{{\operatorname{Pt}}}
\def\today{\ifcase\month\or
 January\or February\or March\or April\or May\or June\or
 July\or August\or September\or October\or November\or December\fi
 \space\number\day, \number\year}
\topmatter
\title  Differential geometry of $\frak g$-manifolds 
\endtitle
\author 
D\. V\. Alekseevsky \\
Peter W\. Michor  
\endauthor
\leftheadtext{\smc D\. V\. Alekseevsky, P\. W. Michor}
\rightheadtext{\smc Differential geometry of $\frak g$-manifolds }
\affil
Erwin Schr\"odinger International Institute of Mathematical Physics, 
Wien, Austria \\
Institut f\"ur Mathematik, Universit\"at Wien, Austria
\endaffil
\address 	 D\. V\. Alekseevsky: 
gen. Antonova 2 kv 99, 117279 Moscow B-279, Russia
\endaddress
\address
P\. W\. Michor: Institut f\"ur Mathematik, Universit\"at Wien,
Strudlhofgasse 4, A-1090 Wien, Austria; and 
Erwin Schr\"odinger International Institute of Mathematical Physics, 
Wien, Austria 
\endaddress
\email Peter.Michor\@esi.ac.at \endemail
\date Received 12 July 1994 \enddate
\thanks Supported by Project P 7724 PHY
of `Fonds zur F\"orderung der wissenschaftlichen Forschung'.
\endthanks
\dedicatory Dedicated to the memory of Professor Franco Tricerri
\enddedicatory
\abstract 
An action of a Lie algebra $\g$ on a manifold $M$ is just a Lie 
algebra homomorphism $\ze:\g\to \X(M)$. 
We define orbits for such an action. In general the space of orbits 
$M/\g$ is not a manifold and even has a bad topology. Nevertheless 
for a $\g$-manifold with equidimensional orbits we treat such 
notions as connection, curvature, covariant differentiation, Bianchi 
identity, parallel transport, basic differential forms, basic 
cohomology, and characteristic classes, which generalize the 
corresponding notions for principal $G$-bundles. As one of the 
applications, we derive a sufficient condition for the projection 
$M\to M/\g$ to be a bundle associated to a principal bundle. 
\endabstract
\endtopmatter


\document

\heading Table of contents \endheading
\noindent 1. Introduction \leaders \hbox to 
     1em{\hss .\hss }\hfill {\eightrm 1}\par
\noindent 2. Lie algebra actions alias ${\frak g}$-manifolds 
     \leaders \hbox to 1em{\hss .\hss }\hfill {\eightrm 4}\par
\noindent 3. Principal connections for Lie algebra actions 
     \leaders \hbox to 1em{\hss .\hss }\hfill {\eightrm 9}\par
\noindent 4. Fr{\accent "7F o}licher-Nijenhuis bracket and curvature 
     \leaders \hbox to 1em{\hss .\hss }\hfill {\eightrm 11}\par
\noindent 5. Homogeneous ${\frak g}$-manifolds \leaders \hbox to 
     1em{\hss .\hss }\hfill {\eightrm 17}\par
\noindent 6. Parallel transport \leaders \hbox to 
     1em{\hss .\hss }\hfill {\eightrm 26}\par
\noindent 7. Characteristic classes for ${\frak g}$-manifolds 
     \leaders \hbox to 1em{\hss .\hss }\hfill {\eightrm 29}\par

\head\totoc\nmb0{1}. Introduction \endhead

Let $\g$ be a finite dimensional Lie algebra and let $M$ be a smooth 
manifold. We say that $\g$ acts on $M$ or that $M$ is a $\g$-manifold 
if there is a Lie algebra homomorphism $\ze=\ze^M:\g\to\X(M)$ from 
$\g$ into the Lie algebra of all vector fields on $M$. 
Many notions and results of the theory of 
$G$-manifolds and of the theory of principal bundles may be extended 
to the category of $\g$-manifolds. This  
is the guideline for our approach to $\g$-manifolds.

Now we describe the structure of the paper and we state some 
principal results.

In section \nmb!{2} notations are fixed and different properties of 
an action of a Lie algebra $\g$ on a manifold $M$ are defined. The 
pseudogroup $\Ga(\g)$ of local transformations generated by an action 
of a Lie algebra $\g$ is considered, and its graph is defined. 
We consider also the groupoid $P(\g)$ of germs of elements from 
$\Ga(\g)$ and under some conditions we may define the adjoint 
representation of $P$ into the adjoint group $\Ad(\g)$ associated 
with the Lie algebra $\g$. Some technical lemmas are proved which 
will be used in section \nmb!{5}. 

In section \nmb!{3} the main definition of a principal connection on 
a $\g$-manifold $M$ is given as such a $\g$-invariant field $\Ph$ of 
endomorphisms of $TM$, whose value in a point $x\in M$ is a projection 
$\Ph_x:T_xM\to \g(x)$ of the tangent space onto the `vertical 
subspace' $\g(x):= \{\ze_X(x):X\in \g\}$. We say, that a principal 
connection $\Ph$ admits a principal connection form if it may be 
represented as $\Ph=\ze_\om$, where $\om$ is a $\g$-valued 
$\g$-equivariant 1-form on $M$ such that 
\roster
\item $T_xM=\g(x)\oplus \ker(\om_x)$ for each $x\in M$,
\item $\g = \g_x\oplus \om(T_xM)$,
\endroster
where $\g_x=\{X\in\g:\ze_X(x)=0\}$ is the isotropy subalgebra. 
Any such form defines a principal connection $\Ph=\ze_\om$. On the 
other hand, a simple example shows that not every principal connection 
admits a principal connection form.

Principal connections may exist only if the action of $\g$ on $M$ has 
constant rank, see proposition \nmb!{3.2} which also gives some 
sufficient conditions for the existence of principal connections.

In order to define the curvature of a principal connection $\Ph$ we 
recall in section \nmb!{4} the definition of the {\it algebraic} bracket 
$[\ph,\ps]^\wedge $ of $\g$-valued differential forms on a manifold 
$M$ which turns the space $\Om(M;\g)$ of such forms into a graded 
Lie algebra. We also recall the definition of the {\it differential} 
Fr\"olicher-Nijenhuis bracket which extends the Lie bracket of vector 
fields to a graded bracket on the space $\Om(M;TM)$ of tangent bundle 
valued differential forms on $M$. An action of a Lie algebra $\g$ on 
$M$ (i\.e\. a homomorphism $\ze:\g\to\X(M)$) induces a linear mapping 
$$
\ze:\Om(M;\g)\to \Om(M;TM).
$$
It is not a homomorphism of graded Lie algebras, but becomes an anti 
homomorphism when it is restricted to the subalgebra 
$\Om_{\text{hor}}^p(M;\g)^\g$ of $\g$-equivariant horizontal forms. 
In general the Fr\"olicher-Nijenhuis bracket $[\ze_\ph,\ze_\ps]$ of two 
$\g$-equivariant forms $\ph,\ps\in\Om(M;\g)$ may be expressed in 
terms of $[\ph,\ps]^\wedge$ and exterior differentials. See 
proposition \nmb!{4.4} for the relevant formulas. 

In section \nmb!{6} we give a local description of a principal 
connection $\Ph$ and of its curvature on a locally trivial 
$\g$-manifold with standard fiber $S$. We show that locally a 
connection is described by a 1-form on the base with values in the 
centralizer $Z_{\X(S)}(\g)$, which may be considered as the Lie 
algebra of infinitesimal automorphisms of the $\g$-manifold $S$. We 
prove that it is isomorphic to the normalizer $N_\g(\g_x)$ of the 
isotropy subalgebra $\g_x$ of a point $x\in S$. As a corollary we 
obtain the existence of a unique principal connection, which is 
moreover flat, under the assumption that $N_\g(\g_x)=0$. 

We treat the case of a homogeneous $\g$-manifold $M$ in section 
\nmb!{5}. First we consider a $\g$-manifold with a free transitive 
action $\ze$ of $\g$ and we remark that the inverse mapping 
$\ka=\ze\i$ is a Maurer-Cartan form (i\.e\. a 1-form that satisfies 
the Maurer-Cartan equation). For a transitive free action of $\g$ on 
a simply connected manifold $M$, we define by means of the graph of 
the pseudogroup a $\g$-equivariant mapping $M\to G$, the `Cartan 
development' of $M$ into the simply connected Lie group with Lie 
algebra $\g$. It is a local diffeomorphism but in general it is 
neither surjective nor injective. As an immediate application, we 
obtain a well defined mapping from any locally flat simply connected 
$G$-structure of finite type into the standard maximally homogeneous 
$G$-structure. It generalizes the developing of a locally flat 
conformal manifold into the conformal sphere.
We prove that on a simply connected $\g$-manifold $M$ with free 
transitive $\g$-action $\ze$ the centralizer of $\ze(\g)$ in the Lie 
algebra $\X(M)$ of all vector fields on $M$ is isomorphic to $\g$. 
The corresponding free transitive action $\hat\ze:\g\to \X(M)$, a Lie 
algebra anti homomorphism, is called the dual action of $\g$. This is 
not true in general, if $M$ is not simply connected.

Let $H\backslash G$ be a homogeneous $G$-manifold for a Lie group $G$ with 
isotropy group $H\subset G$ of a point $o$. Then $p:G\to H\backslash G$ 
is a principal $H$-bundle with left principal $H$-action, and $G$ 
acts from the right by automorphisms of principal bundles. 
Let $\ka\in\Om^1(G,\g)$ be the right invariant Maurer-Cartan form on 
$G$, associated with the left action of $G$ on itself.
Then any reductive decomposition $\g=\h\oplus \frak m$ with $\Ad(H)\frak m= 
\frak m$ defines a $G$-invariant principal connection 
$\om:=\operatorname{pr}_\h \o \ka$ for the principal bundle 
$p:G\to H\backslash G$, and a $G$-invariant displacement form $\th:= 
\operatorname{pr}_{\frak m}\o \ka$. Any $G$-invariant principal 
connection of $p:G\to H\backslash G$ has this form.

We generalize these classical results in section \nmb!{5} to the case 
of a homogeneous $\g$-manifold $M$. The role of the principal bundle 
$p:G\to H\backslash G$ is taken by the manifold $P$ of germs of 
transformations of the pseudogroup $\Ga(\g)$, at a fixed point. We 
prove that the principal connection forms on the $\g$-manifold $M$ 
correspond bijectively to the invariant principal connections of the 
principal bundle $P\to M$. 

For a locally trivial $\g$-manifold $M$ with a principal connection 
$\Ph$ we define the horizontal lift of vector fields on the orbit 
space $N$ and the parallel transport along a smooth curve on $N$. The 
parallel transport however is only locally defined. If the parallel 
transport is defined on the whole fiber along any smooth curve, then 
the connection is called complete. We show that any principal 
connection is complete if all vector fields in the centralizer 
$Z_{\X(S)}(\g)$ are complete.

As final result in this section we prove the following: If a locally 
trivial $\g$-manifold $M$ with standard fiber $S$ admits a complete 
principal connection $\Ph$, whose holonomy Lie algebra consists of 
complete vector fields on $S$, then the bundle $M\to N=M/\g$ is 
isomorphic to the bundle $P[S]=P\x_H S$ associated to a principal 
$H$-bundle $P\to N$, where $H$ is the holonomy group. Moreover, the 
connection $\Ph$ is induced by a principal connection on $P$.

In the last section \nmb!{7} we assume that the $\g$-manifold $M$ 
admits not only a principal connection $\Ph$, but also a principal 
connection form $\om\in\Om^1(M;\g)$ with curvature form $\Om$. We 
define the Chern-Weil homomorphism $\ga$ from the algebra 
$S(\g^*)^\g$ of $\ad(\g)$-invariant polynomials on $\g$ into the 
algebra $\Om^{\text{closed}}(M)^\g$ of $\g$-invariant closed forms on 
$M$. We prove that for any $f\in S(\g^*)^\g$ the cohomology class 
$[\ga(f)]$ depends only on $f$ and the $\g$-action. If the action of 
$\g$ is free the image of $\ga$ consists of horizontal forms. The 
associated cohomology classes are basic and may be considered as 
characteristic classes of the $\g$-manifold, or of the `bundle below 
$M$', even if the action of $\g$ is not locally trivial.
If on the other hand $M$ is a homogeneous $\g$-manifold, our 
cohomology classes are characteristic classes of the `bundle above 
$M$', the principal bundle $P\to M$ consisting of germs of 
pseudogroup transformations constructed in  \nmb!{5.8}.

\head\totoc\nmb0{2}. Lie algebra actions alias $\g$-manifolds \endhead

\subhead\nmb.{2.1}. Actions of Lie algebras on a manifold  \endsubhead
Let $\g$ be a finite dimensional Lie algebra and let $M$ be a smooth 
manifold. We say that \idx{\it $\g$ acts on $M$} or that $M$ is a 
\idx{\it $\g$-manifold} if there is a Lie algebra homomorphism 
$\ze:\g\to \X(M)$, from $\g$ into the  Lie algebra of all vector 
fields on $M$. 

If we have a right action of a Lie group on $M$, then the 
fundamental vector field mapping is an action of the corresponding  
Lie algebra on $M$.

\proclaim{Lemma} 
If a Lie algebra $\g$ acts on a 
manifold $M$, then it spans an integrable distribution 
on $M$, which need not be of constant rank. So through each point of 
$M$ there is a unique maximal leaf of that distribution; we also call 
it the $\frak g$-orbit through that point. It is an initial 
submanifold of $M$ in the sense that a mapping from a manifold into 
the orbit is smooth if and only if it is smooth into $M$, see 
\cit!{9},~2.14ff.
\endproclaim

\demo{Proof} See \cit!{19} or \cit!{20} for integrable distributions of 
non-constant rank, or \cit!{9}, 3.25. 
Let $\Fl^\xi_t$ denote the flow of a vector field $\xi$. One may 
check easily that for $X,Y\in \g$ we have
$\tfrac{d}{dt}(\Fl^{\ze_X}_{-t})\,\ze(e^{t\ad(X)}Y)=0$,
which implies 
$(\Fl^{\ze_X}_t)^*\ze_Y = \ze(e^{t\ad(X)}Y)$.
So condition (2) 
of theorem \cit!{9},~3.25 is satisfied and all assertions follow.
\qed\enddemo

An action of a Lie algebra $\g$ on a manifold $M$ may have the 
following properties:
\roster
\item It is called \idx{\it effective} if $\ze:\g\to \X(M)$ is 
       injective. So for each  $X\in\g$ there is some $x\in M$ such 
       that $\ze_X(x)\ne0$.
\item The action is called \idx{\it free} if for each $x\in M$ the 
       mapping $X\mapsto \ze_X(x)$ is injective. Then the 
       distribution spanned by $\ze(\g)$ is of constant rank.
\item The action is called \idx{\it transitive} if for each $x\in M$ 
       the mapping $X\mapsto \ze_X(x)$ is surjective onto $T_xM$. If 
       $M$ is connected then it is the only orbit and 
       we call $M$ a \idx{\it homogeneous $\g$-space}.
\item The action is called \idx{\it complete} if each fundamental 
       vector field $\ze_X$ is complete, i\.e\. it generates a global flow.
       In this case the action can be integrated to a right action of 
       a connected Lie group $G$ with Lie algebra $\g$, by a result 
       of Palais, \cit!{15}.
\item The action is said to be of \idx{\it constant rank} $k$ if all 
       orbits have the same dimension $k$.
\item We call it an \idx{\it isostabilizer action} if all the 
       isotropy algebras $\g_x:= \ker(\ze_x:\g\to T_xM)$ are 
       conjugate in $\g$ under the connected adjoint group. An 
       isostabilizer action is of constant rank.
\item The action is called \idx{\it locally trivial} if there exists a 
       connected manifold 
       $S$ with a transitive action of $\g$ on $S$, a submersion 
       $p:M\to N$ onto a smooth manifold with trivial 
       $\g$-action, such that for each point $x\in N$ there exists an 
       open neighborhood $U$ and a $\g$-equivariant diffeomorphism 
       $\ph:p\i(U)\to U\x S$ with $p\o \ph=p$. By 
       `$\g$-equivariant' we mean that for each $X\in\g$ the 
       fundamental vector fields $\ze^M_X|p\i(U)$ and $0\x\ze^S_X$ 
       are $\ph$-related: 
       $T\ph\o \ze^M_X|p\i(U)=(0\x\ze^S_X)\o\ph$. 
       A pair like $(U,\ph)$ is called a \idx{\it bundle chart}. Note 
       that $N$ is canonically isomorphic to the space of orbits 
       $M/\g$, and that a locally trivial action is isostabilizer and 
       of constant rank.
\item A free and locally trivial $\g$-action is called a 
       \idx{\it principal action}. A complete principal $\g$-action 
       can be integrated to an almost free action of a Lie group $G$ 
       (i\.e\. with discrete isotropy groups). If its action is free, 
       it defines a principal bundle $p:M\to M/G$, and the action 
       of $\g$ on $M$ is the associated action of the Lie algebra of 
       $G$. This explains the name.
\endroster

In the general case, we will consider a locally trivial 
$\g$-manifold $M$ as some generalization of the notion of a 
principal $G$-bundle, and we will extend to this case some of 
the main differential geometric constructions of the geometry 
of principal bundles.

A smooth mapping $f:M\to N$ between $\g$-manifolds $M$ and $N$ is 
called \idx{\it $\g$-equivariant} if for each $X\in \g$ the 
fundamental vector fields $\ze^M_X$ and $\ze^N_X$ are $f$-related:
$Tf\o \ze^M_X = \ze^N_X\o f$. In view of \cit!{9}, section 47, we may 
also say, that the \idx{\it generalized Lie derivative} of $f$ is 
zero: 
$$
\L_Xf=\L_{\ze^M_X,\ze^N_X}f:= \ze^N_X\o f - Tf\o \ze^M_X = 0
$$

Note that the integrable distribution of a $\g$-manifold $M$ of 
constant rank is a special case (in a certain sense the simplest 
case) of a foliation. To make this statement more precise we define 
the \idx{\it degree of cohomogeneity} of a foliation (integrable 
distribution) $\Cal D$ on $M$ as the minimum of the difference 
between the rank of $\Cal D$ and the rank of a $\g$-manifold 
structure on $M$ where $\g$ runs through all finite dimensional 
subalgebras of constant rank  in the Lie algebra $\X(\Cal D)$ of 
global vector fields on $M$ which are tangent to $\Cal D$:
$$
\min\{\operatorname{rank}(\Cal D)-\operatorname{rank}(\g): 
\g\subset \X(\Cal D), \dim(\g)<\infty\}
$$
Then we may say that the foliation associated with a $\g$-manifold of 
constant rank has degree of cohomogeneity 0, or that it is a 
`homogeneous foliation'.

\subhead\nmb.{2.2}. The pseudogroup of a $\g$-manifold \endsubhead
Let $M$ be a $\g$-manifold which we assume to be effective and 
connected. Local flows of fundamental vector fields, restricted to 
open subsets, and their compositions, form the \idx{\it pseudogroup} 
$\Ga(\g)$ of the $\g$-action. 

Let us first recall the following definition: A pseudogroup of 
diffeomorphisms of the manifold $M$ is a set $\Ga$ consisting of 
diffeomorphisms $\ph:U\to V$ between connected open subsets of $M$, 
subject to the following conditions:
\roster
\item If $\ph:U\to V$ is an element of $\Ga$ then also 
       $\ph\i:V\to U$.
\item If $\ph:U\to V$ and $\ps:V\to W$ are elements of $\Ga$ then 
       also the composition $\ps\o\ph:U\to W$ is in $\Ga$.
\item If $\ph:U\to V$ is an element of $\Ga$ then also its 
       restriction to any connected open subset $U_1\subset U$ is an 
       element of $\Ga$.
\item If $\ph:U\to V$ is a diffeomorphism between connected open 
       subsets of $M$ which coincides on an open neighborhood of each 
       of its points with an element of $\Ga$ then also $\ph$ is in 
       $\Ga$.
\endroster
Now in more details $\Ga(\g)$ consists of diffeomorphisms of the 
following form:
$$
\Fl^{\ze_{X_n}}_{t_n}\o\dots
     \o\Fl^{\ze_{X_2}}_{t_2}\o\Fl^{\ze_{X_1}}_{t_1}|U
\tag5$$
where $X_i\in\g$, $t_i\in\Bbb R$, and $U\subset M$ are such that 
$\Fl^{\ze_{X_1}}_{t_1}$ is defined on $U$, $\Fl^{\ze_{X_2}}_{t_2}$ is 
defined on $\Fl^{\ze_{X_1}}_{t_1}(U)$, and so on.

\subhead\nmb.{2.3}. The graph of the pseudogroup of a $\g$-manifold 
\endsubhead
Let $M$ be a connected $\g$-manifold. Let $G$ be a connected Lie group 
with Lie algebra $\g$. We consider the  distribution of rank $\dim\g$ 
on $G\x M$ which is given by
$$
\{(L_X(g),\ze^M_X(x)):(g,x)\in G\x M, X\in\g\}\subset TG\x TM,\tag1
$$
where $L_X$ is the left invariant vector field on $G$ generated by 
$X\in\g$. Obviously this distribution is integrable and thus we may 
consider the foliation induced by it, which we will call the 
\idx{\it graph of the pseudogroup $\Ga(\g)$}.
Note that the flow of the vector field $(L_X,\ze^M_X)$ on $G\x M$ is 
given by 
$$
\Fl^{(L_X,\ze^M_X)}_t(x,g)=(g.\exp^G(tX),\Fl^{\ze_X}_t(x)).
$$ 
In the sense of \cit!{9},~section~9, 
this foliation is the horizontal foliation for a flat connection 
of the trivial fiber bundle $G\x M\to G$.
The first projection $\operatorname{pr}_1:G\x M\to G$, when 
restricted to a leaf, is locally a diffeomorphism. For $x\in M$ we 
consider the leaf $L(x)$ through $(e,x)\in G\x M$. Then 
$W_x:=\operatorname{pr}_1(L(x))$ is a connected open set in $G$.

In particular we may use the theory of 
parallel transport \cit!{9},~9.8:
Let $c:(a,b)\to G$ be a piecewise smooth curve with $0\in(a,b)$ and $c(0)=g$.
Then there is an open subset $V$ of $\{g\}\x M\x \{0\}$ in 
$\{g\}\x M\x \Bbb R$ and a smooth mapping $\Pt_c:V\to G\x M$ such that:
\roster
\item $\operatorname{pr}_1(\Pt(c,(g,x),t)) = c(t)$ if defined, and
     $\Pt(c,(g,x),0) = (g,x)$.
\item $\frac d{dt}\Pt(c,(g,x),t)$ is tangent to the graph foliation.
\item Reparametrization invariance: If $f: (a',b') \to   (a,b)$ is 
       piecewise smooth with  $0 \in (a',b')$, then $\Pt(c,(g,x),f(t)) 
       = \Pt(c \o f, \Pt(c,u_x,f(0)),t)$ if defined.    
\item $V$ is maximal for properties \therosteritem1 and \therosteritem2.
\item If the curve $c$ depends smoothly on further parameters then 
      $\Pt(c,(g,x),t)$ depends also smoothly on those parameters.
\endroster
Now let $c:[0,1]\to G$ be piecewise smooth with $c(0)=e$, and let us 
assume that for some $x\in M$ the parallel transport $\Pt(c,(e,x),t)$ 
is defined for all $t\in[0,1]$. Then in particular 
$c([0,1])\subset W_x$. Since $\{(e,x)\}\x[0,1]\subset V$ the parallel 
transport $\Pt(c,1)$ is defined on an open subset $\{e\}\x U$ of 
$(e,x)$, and by \therosteritem3 it is a diffeomorphism onto its 
image $\{c(1)\}\x U'$. We may choose $U$ maximal 
with respect to this property. 
Since the connection is flat the parallel transport depends on the 
curve $c$ only up to small (liftable) homotopies fixing end points, since 
$\Pt(c,(e,x))$ is just the unique lift over the local diffeomorphism 
$\operatorname{pr}_1:L(x)\to W_x$. We put
$$
\ga_x(c):= \operatorname{pr}_2\o\Pt(c,1)\o\operatorname{ins}_e:
     U\to \{e\}\x U\to \{c(1)\}\x U'\to U',
$$
so $\ga_x(c)$ is the parallel transport along $c$, from the fiber 
over $e$ to the fiber over $c(1)$, viewed as a local diffeomorphism 
in $M$.
Since $c$ is homotopic within $W_x$ to a finite sequence of left  
translates of 1-parameter subgroups, this parallel transport is a 
composition of a sequence of flows of fundamental vector fields,
so $\ga_x(c)$ is an element of the pseudogroup $\Ga(\g)$ on $M$.
So $\ga_x$ is a mapping from the set of homotopy classes fixing 
end points of curves starting at $e$ in $W_x$ into the pseudogroup 
$\Ga(\g)$. 

Conversely each element of $\Ga(\g)$ of the form \nmb!{2.2}.\thetag5 
applied to $x\in U$ is the parallel transport of $(e,x)$ along the 
corresponding polygonial arc consisting of left translates of 
1-parameter subgroups: first $[0,t_1]\ni t\mapsto \exp(tX_1)$, then 
$[t_1,t_2]\ni t\mapsto \exp(t_1X_1).\exp((t-t_1)X_2)$, and so on.
Thus we have proved:

\proclaim{Lemma} Let $M$ be a connected $\g$-manifold. 
Then any element $\ph:U\to V$ of the pseudogroup $\Ga(\g)$ is of the 
form $\ph=\ga_x(c)$ for $x\in U$ and a smooth curve $c:[0,1]\to W_x$. 
\endproclaim 

\proclaim{\nmb.{2.4}. Lemma} Let $M$ be a $\g$-manifold. Assume that 
for a point $x\in M$ the Lie algebra homomorphism
$\operatorname{germ}_x\o\ze:\g\to \X(M)_{\text{germs at }x}$ is 
injective. Then for each $\ph\in\Ga(\g)$ which is defined near $x$ 
there is a unique automorphism $\Ad(\ph\i):\g\to\g$ satisfying 
$\ph^*\ze_X=\ze_{\Ad(\ph\i)X}$ for all $X\in \g$. 
\endproclaim

This mapping $\Ad$ generalizes the adjoint representation of a Lie 
group.

\demo{Proof}
One may check easily that for $X,Y\in \g$ we have
$\tfrac{d}{dt}(\Fl^{\ze_X}_{-t})^*\,\ze(e^{t\ad(X)}Y)=0$ for all $t$ 
for which the flow is defined.
This implies
$$
(\Fl^{\ze_X}_t)^*\ze_Y = \ze(e^{t\ad(X)}Y).
\tag1$$
We may apply \thetag1 iteratively to elements of $\Ga(\g)$ of the 
form \nmb!{2.2}.\thetag5 and thus we get
$$
\ga_x(c)^*\ze_Y = \ze_{\Ad(c(1))Y}
\tag2$$
for each smooth curve starting from $e$ in $W_x$ which is liftable to 
$L(x)$ in the setting of \nmb!{2.3}. By the assumption, equation 
\thetag2 now implies that $\Ad(c(1))$ depends only on 
$\ga_x(c)\in\Ga(\g)$ and we call it $\Ad(\ga_x(c)\i)$. We use the 
inverse so that $\Ad$ becomes a `homomorphism' in \nmb!{2.5} below.
\qed\enddemo

\subhead\nmb.{2.5}. Adjoint representation  \endsubhead
Let $M$ be a $\g$-manifold with pseudogroup $\Ga(\g)$ such that for 
each $x\in M$ the homomorphism 
$\operatorname{germ}_x\o \ze:\g \to \X(M)_{\text{germs at }x}$ is 
injective. We denote by 
$P_x(\g)$ the set of all germs at $x\in M$ of transformations in 
$\Ga(\g)$ which are defined at $x$. Then the set $P(\g):= 
\bigcup_{x\in M}P_x(\g)$ with the obvious partial composition is a 
groupoid. By lemma \nmb!{2.4} we have a well defined representation
$$
\Ad:P \to \Ad(\g)
$$ 
with values in the adjoint group,
and we call it the \idx{\it adjoint representation} of the groupoid $P$.

\proclaim{\nmb.{2.7}. Lemma}                                      
Let $M$ be a $\g$-manifold. Let $c:[0,1]\to M$ be a smooth curve in 
$M$ with values in one $\g$-orbit. Then there exists a smooth mapping
$\ph:[0,1]\x M\supset U\to M$ such that $\ph_t\in\Ga(\g)$ for each 
$t$, $[0,1]\x\{c(0)\}\subset U$, and $c(t)=\ph(t,c(0))$ for all $t$.
\endproclaim
Since each $\g$-orbit is an initial submanifold we may equivalently 
assume that $c$ is a smooth curve in a $\g$-orbit.

\demo{Proof} Let us call $c(0)=x$.
Since $c'(t)\in\ze_{c(t)}(\g)$ it is easy to get a smooth curve 
$b:[0,1]\to \g$ such that $c'(t)=\ze_{b(t)}(c(t))$. We may choose for 
example $b(t)=(\ze_{c(t)}|\ker(\ze_{c(t)})^\bot)\i(c'(t))$ with 
respect to any inner product on $\g$.
Let $G$ be a Lie 
group with Lie algebra $\g$, let $g(t)$ be the integral curve of the 
time dependent vector field $(t,g)\mapsto L_{b(t)}(g)$ with $g(0)=e$.
Then $(g(t),c(t))$ is a smooth curve in $G\x M$ which is tangent to 
the graph foliation of the pseudogroup $\Ga(\g)$ and thus it lies in 
the leaf through $(e,x=c(0))$. From \nmb!{2.3} we see that 
$\ga_x(g|[0,t])=\ph_t\in\Ga(\g)$, where $\ph$ is the evolution 
operator of the time dependent vector field 
$(t,x)\to \ze_{(b(t))}(x)$ on $M$.
\qed\enddemo

\subhead\nmb.{2.6}. Isotropy groups \endsubhead Let $M$ be a 
connected $\g$-manifold and let $x\in M$. Let us denote by 
$\Ga(\g)_{x}$ the group of all germs at $x$ of elements of the 
pseudogroup $\Ga(\g)$ fixing $x$. It is called the 
\idx{\it isotropy group}. Its natural representation into the 
space $J^k_{x}(M,\Bbb R)_0$ of $k$-jets at $x$ of functions 
vanishing at $x$ is called the \idx{\it isotropy representation of 
order} $k$. In general the isotropy representation of any order may 
have a nontrivial kernel. The simplest example is provided by the Lie 
algebra action defined by one vector field which is flat (vanishes 
together with all derivatives) at $x$. We remark that this cannot 
happen if $\Ga(\g)$ is a `Lie pseudogroup' (defined by a system of 
differential equations); in particular if $\ze(\g)$ is the algebra of 
infinitesimal automorphisms of some geometrical structure.

Consider the 
following diagram
$$\CD
\g   @>{\ze}>>   \X(M) \\
@V{\ze_{x}}VV  @VV{\text{germ}_{x}}V\\
T_{x}M @<<{\ev_{x}}<  \X(M)_{\text{germs at }x} 
\endCD$$
The kernel of the linear mapping $\ze_{x}:\g\to T_{x}M$ is 
denoted by $\g_{x}$ and it is called the \idx{\it isotropy algebra} 
at $x$. The kernel of the Lie algebra homomorphism 
$\operatorname{germ}_{x}\o \ze:\g \to \X(M)_{\text{germs at }x _0}$ 
is denoted by $\g_{\text{germ}_{x}=0}$; it is an ideal of $\g$ 
contained in the isotropy algebra $\g_{x}$.

\proclaim{Lemma}
In this setting, $\Ga(\g)_{x}$ is a Lie group (not necessarily 
second countable) whose Lie algebra is anti isomorphic to 
the quotient $\g_{x}/\g_{\text{germ}_{x}=0}$ of the isotropy 
algebra $\g_{x}$.

If the Lie algebra homomorphism 
$\operatorname{germ}_x\o\ze:\g\to \X(M)_{\text{germs at }x}$ 
is injective, then there is a canonical representation 
$\Ad:\Ga(\g)_{x}\to \operatorname{Aut}(\g)$ which leaves invariant 
the isotropy subalgebra $\g_{x}$ and coincides on $\g_x$ with the 
adjoint representation of $\Ga(\g)_x$.
\endproclaim

\demo{Proof}
As in \nmb!{2.3} we consider again the graph foliation of the 
$\g$-manifold $M$ on $G\x M$, where $G$ is a connected Lie group with 
Lie algebra $\g$, the leaf $L(x)$ through $(e,x)$ of it, and the open 
set $W_{x}=\operatorname{pr}_1(L(x))\subset G$. Let $G_{x}$ be the 
connected subgroup of $G$ corresponding to the isotropy algebra 
$\g_{x}$. Then $G_{x}$ is contained in $W_{x}$ since for a smooth 
curve $c:[0,1]\to G_{x}$ the curve $(c(t),x)$ in $G\x M$ is tangent 
to the graph foliation; each curve in $G_x$ and even each homotopy in 
$G_x$ is liftable to $L(x)$. The universal cover of $G_x$ may be 
viewed as the space of homotopy classes with fixed ends, of smooth 
curves in $G_x$ starting from $e$. So by assigning the germ at $x$ of 
$\ga_{x}(c)\in \Ga(\g)_{x}$ to the homotopy class of a curve $c$ in 
$G_x$ starting from $e$, we get a group homomorphism from the 
universal cover of $G_{x}$ into $\Ga(\g)_{x}$. Its tangent mapping at 
the identity is $-\operatorname{germ}_{x}\o\ze|\g_{x}$.
Let us denote by $\Ga(\g)_{x}^0$ the image of this group 
homomorphism.

Let $\ph_t$ be a smooth curve in the group of germs $\Ga(\g)_{x}$, 
with $\ph_0=Id$ in the sense that $(t,x)\to \ph_t(x)$ is a smooth 
germ. Then $(\tfrac d{dt}\ph_t)\o\ph_t\i$ is the germ at $x$ of a 
time dependent vector field with values in the distribution $\g(M)$ 
spanned by $\g$, which vanishes at $x$. So it has values in the set of 
germs at $x$ of $\ze_{\g_{x}}$, and thus $\ph_t$ is in 
$\Ga(\g)^0_{x}$, see the proof of \nmb!{2.7}. So  the normal subgroup 
of $\Ga(\g)_x$ of those elements which may be connected with the 
identity by a smooth curve in  $\Ga(\g)_x$, coincides with 
$\Ga(\g)_x^0$, and the latter is a normal subgroup. 

If we declare the Lie group $\Ga(\g)_{x}^0$ to be open in 
$\Ga(\g)_{x}$ we get a Lie group structure on $\Ga(\g)_{x}$.

The statement about the adjoint representation follows immediately 
from lemma \nmb!{2.4}.
\qed\enddemo

The anti isomorphism in this lemma comes because $\Ga(\g)_x$ acts 
from the left on $M$, so the fundamental vector field mapping of this 
action should be a Lie algebra anti isomorphism, see \cit!{9},~5.12. 
Since we started from a Lie algebra homomorphism $\ze:\g\to\X(M)$, 
the pseudogroup should really act from the right; so it should be 
viewed as an 
abstract pseudogroup and not one of transformations. We decided not 
to do this, but this will cause complicated sign conventions, 
especially in theorem \nmb!{5.8} below.

\head\totoc\nmb0{3}. Principal connections for Lie algebra actions \endhead

\subhead\nmb.{3.1}. Principal connections \endsubhead
Let $M$ be a $\g$-manifold. A vector valued 1-form 
$\Ph\in\Om^1(M;TM)$, i\.e\. a vector bundle homomorphism 
$\Ph:TM\to TM$, is called a \idx{\it connection} for the $\g$-action 
if for each $x\in M$ the mapping $\Ph_x:T_xM\to T_xM$ is a projection 
onto $\g(x)=\ze(\g)(x)\subset T_xM$. The connection is called 
\idx{\it principal} if it is $\g$-equivariant, i\.e\. 
if for each $X\in \g$ the Lie derivative 
vanishes: $\L_{\ze_X}\Ph=[\ze_X,\Ph]=0$, where $[\quad,\quad]$ is the 
Fr\"olicher-Nijenhuis bracket, see \nmb!{4.4} below. The distribution 
$\ker(\Ph)$ is called the \idx{\it horizontal distribution} of the 
connection $\Ph$.

A Lie algebra valued 1-form $\om\in\Om^1(M;\g)$ is called a 
\idx{\it principal connection form} if the following conditions 
are satisfied:
\roster
\item $\om$ is $\g$-equivariant, i\.e\. for all $X\in\g$ we have 
       $\L_{\ze_X}\om=-\ad(X)\o \om$.
\item For any $x\in M$ we have $\ze_x=\ze_x\o\om_x\o\ze_x: 
       \g\to T_xM\to\g\to T_xM$.
\endroster
Thus for any $x\in M$ the kernel $\ker(\om_x)$ is a complementary 
subspace to the vertical space $\g(x)$ and 
the mapping $\om_x:\g(x)\cong \g/g_x\to \g$ is 
a right inverse to the projection $\g\to \g/\g_x$.             

Any principal connection form $\om$ defines a principal connection 
$\Ph:= \ze_\om$. The converse statement is not true in general as  
example \nmb!{5.7} shows.

\proclaim{\nmb.{3.2}. Proposition} 
Let $M$ be a $\g$-manifold.

1. If $M$ admits a principal connection then the action of $\g$ on 
$M$ has constant rank.

2. Let us assume conversely that the $\g$-action is of constant rank. 
Then $M$ admits a principal connection if any of the following 
conditions is satisfied:
\roster
\item The action is locally trivial.
\item There exists a $\g$-invariant Riemannian metric on $M$.
\item The $\g$-action is induced by a proper 
       action of a Lie group $G$ with Lie algebra $\g$.
\endroster

3. Let the $\g$-action be locally trivial with standard fiber $S$, 
a homogeneous $\g$-space which admits a principal connection form.
Then $M$ admits even a principal connection form $\om$.
\endproclaim

For assertion 3, see example \nmb!{5.7} for conditions assuring the 
existence of principal connection forms on the standard fiber $S$:
the isotropy subalgebra $\g_x$ of some point $x\in S$ admits an 
$\ad(g_x)$-invariant complement $\frak m$ in $\g$ which is also 
invariant under the isotropy representation of the pseudogroup 
$\Ga(\g)$ generated by $\g$.

\demo{Proof}
1. If a principal connection $\Ph$ exists, it is a projection onto the 
distribution spanned by $\g$ (which we will call the vertical 
distribution sometimes), and its rank cannot fall locally. But the 
rank of the complementary projection $\ch:=Id_{TM}-\Ph$ onto the 
kernel of $\Ph$ also cannot fall locally, so the vertical 
distribution $\g(M)$ is locally of constant rank.

2. First of all, we have the implications \therosteritem3 $\implies$ 
\therosteritem2 $\implies$ \therosteritem1. The first implication is a 
theorem of Palais \cit!{16}. The implication \therosteritem2 
$\implies$ \therosteritem1 may be proved as for an action of a Lie 
group that preserves a Riemannian metric, using a slice. 
Hence, we may 
assume that $M$ is a connected locally trivial $\g$-manifold.

Let $(U_\al,\ph_\al:p\i(U_\al)\to U_\al\x S)$ be a family of 
principal charts such that $(U_\al)$ is an open cover of $M$. Put
$\Ph_\al(T\ph_\al\i(\xi_x,\et_s))=T\ph_\al\i(0_x,\et_s)$ for 
$\xi_x\in T_xU_\al$ and $\et_s\in T_sS$.
Obviously that $\Ph_\al$ is a principal 
connection on $p\i(U_\al)$. Now let $f_\al$ be a smooth 
partition of unity on $N$ which is subordinated to the open cover 
$(U_\al)$. Then $\Ph:= \sum_\al (f_\al\o p)\Ph_\al$ is a principal 
connection on $M$.

3. This is proved similarly as 2, starting from a principal connection 
form on the standard fiber $S$.
\qed\enddemo

\subhead\nmb.{3.3} \endsubhead
Let $M$ be a $\g$-manifold. If a principal connection $\Ph$ exists 
then the distribution $\g(M)$ spanned by $\g$ is of constant rank and 
thus a vector bundle over $M$, and $\Ph$ factors to a 
$\g$-equivariant right inverse
of the vector bundle epimorphims $TM \to TM/\g(M)$.
Let us consider the following sequence of families of vector bundles
over $M$, where
$\operatorname{iso}(M):= \bigcup_{x\in M}\{x\}\x \g_x = \ker(\ze^M)$ 
is the \idx{\it isotropy algebra bundle} over $M$:
$$
\operatorname{iso}(M) @>>> M\x \g @>{\ze^M}>> TM @>>> TM/\g(M)
$$ 
Then a principal connection form $\om$ induces a $\g$-equivariant 
right inverse on it's image of the vector bundle homomorphism 
$\ze^M:M\x \g\to TM$, so it satisfies $\ze\o\om\o\ze=\ze$.

\head\totoc\nmb0{4}. Fr\"olicher-Nijenhuis bracket and curvature \endhead

\subheading{\nmb.{4.1} Products of differential forms}
Let $\rh:\g\to \frak g\frak l(V)$ be a representation of a Lie 
algebra $\g$ in a finite dimensional 
vector space $V$ and let $M$ be a smooth manifold.

For $\ph\in\Om^p(M;\g)$ and $\Ps\in\Om^q(M;V)$ we define the 
form	$\rh^\wedge (\ph)\Ps\in\Om^{p+q}(M;V)$	by 
$$\multline
(\rh^\wedge (\ph)\Ps)(\row X1{p+q}) =\\
= \frac 1{p!\,q!} \sum_{\si} \text{sign}(\si)
 	\rh(\ph(\row X{\si1}{\si p}))\Ps(\row X{\si(p+1)}{\si(p+q)}).
\endmultline$$
Then $\rh^\wedge (\ph):\Om^*(M;V)\to \Om^{*+p}(M;V)$ is a graded 
$\Om(M)$-module homomorphism of degree $p$.

Recall also that $\Om(M;\g)$ is a graded Lie algebra with the bracket 
$[\quad,\quad]^\wedge = [\quad,\quad]^\wedge_\g$ given by
$$\multline
[\ph,\ps]^\wedge(\row X1{p+q}) = \\
= \frac 1{p!\,q!} \sum_{\si} \text{sign}\si\,
	[\ph(\row	X{\si1}{\si p}),\ps(\row X{\si(p+1)}{\si(p+q)})]_{\g},
\endmultline$$
where $[\quad,\quad]_{\g}$ is the bracket in $\g$.
One may easily check that for the graded commutator in 
$\End(\Om(M;V))$ we have
$$
\rh^\wedge ([\ph,\ps]^\wedge ) = 
     [\rh^\wedge (\ph),\rh^\wedge (\ps)] = 
     \rh^\wedge (\ph)\o \rh^\wedge (\ps) - (-1)^{pq} 
     \rh^\wedge (\ps)\o \rh^\wedge (\ph)
$$
so that $\rh^\wedge :\Om^*(M;\g) \to \End^*(\Om(M;V))$ is a 
homomorphism of graded Lie algebras.

For any vector space $V$ let
$\bigotimes V$ be the tensor algebra generated by $V$.
For $\Ph,\Ps\in \Om(M;\bigotimes V)$ we will use the associative 
bigraded product
$$\multline
(\Ph\otimes_\wedge \Ps)(\row X1{p+q}) = \\
= \frac 1{p!\,q!} \sum_{\si} \text{sign}(\si)
 	\Ph(\row X{\si1}{\si p})\otimes \Ps(\row X{\si(p+1)}{\si(p+q)})
\endmultline$$

\subhead\nmb.{4.2}. Basic differential forms \endsubhead
Let $M$ be a $\g$-manifold. 
A differential form $\ph\in\Om^p(M;V)$ with values in a vector space 
$V$ (or even in a vector bundle over $M$) is called \idx{\it horizontal}
if it kills all fundamental vector fields $\ze_X$, i\.e\. if 
$i_{\ze_X}\ph = 0$ for each $X\in \g$.

If moreover $\rh:\g\to \g\frak l(V)$ is a representation of the Lie 
algebra $\g$ in $V$, then differential form $\ph\in\Om^p(M;V)$ is 
called \idx{\it $\g$-equivariant} if for the Lie derivative along 
fundamental vector fields we have:
$\L_{\ze_X}\ph = -\rh(X)\o \ph$ for all $X\in\g$.

Let us denote by $\Om_{\text{hor}}^p(M;V)^\g$ the space of all 
$V$-valued differential forms on $M$ which are horizontal and 
$\g$-equivariant. It is called the space of \idx{\it basic $V$-valued 
differential forms} on  the $\g$-manifold $M$. If the $\g$-manifold 
$M$ has constant rank and the action of $\g$ defines a foliation, 
scalar valued basic forms are the usual 
basic differential forms of the foliation, see 
e\.g\. \cit!{13}.

Note that the graded Lie module structure $\rh^\wedge$ from 
\nmb!{4.1} restricts to a graded Lie module structure 
$\rh^\wedge :\Om(M;\g)^\g\x \Om(M;V)^\g \to \Om(M;V)^\g$.
It is also compatible with the requirement of horizontality.

The exterior differential $d$ acts on 
$(\Om(M;\g),[\quad,\quad]^\wedge)$ as a graded derivation of degree~1.
It preserves the subalgebra $\Om(M;\g)^\g$ of $\g$-invariant forms, 
but it does not preserve the subalgebra $\Om_{\text{hor}}^p(M;\g)^\g$ 
of $\g$-valued basic forms.

For $\ph\in\Om^p(M;\g)$ we consider the tangent bundle 
valued differential form $\ze_\ph\in \Om^p(M;TM)$ which is given for 
$\xi_i\in T_xM$ by 
$$
(\ze_\ph)_x(\xi_1,\dots,\xi_p):= \ze_{\ph_x(\xi_1,\dots,\xi_p)}(x).
$$

\subhead\nmb.{4.3}. Fr\"olicher-Nijenhuis bracket \endsubhead
Let $M$ be a smooth manifold.
We shall use now the Fr\"olicher-Nijenhuis bracket 
$$
[\quad,\quad]:\Om^p(M;TM)\x\Om^q(M;TM)\to \Om^{p+q}(M;TM)
$$
as guiding line for the further developments, since it is a natural 
and convenient way towards connections, curvature, and Bianchi 
identity, in many settings. See \cit!{9}, sections~8--11, as a convenient 
reference for this. We repeat here the global formula for the 
Fr\"olicher-Nijenhuis bracket from \cit!{9},~8.9: 
For $K\in \Om^k(M;TM)$ and $L\in \Om^\ell(M;TM)$ 
we have for the Fr\"olicher-Nijenhuis bracket
$[K,L]$ the following formula, where the $\xi_i$ are vector fields
on $M$.
$$\align
[K,L]&(\row \xi1{k+\ell}) = \tag1\\
&= {\tsize\frac 1{k!\,\ell!}}\sum_{\si} \sign\si\;
     [K(\row \xi{\si1}{\si k}),L(\row \xi{\si (k+1)}{\si(k+\ell)})]\\
& + {\tsize\frac{-1}{k!\,(\ell-1)!}}\sum_\si \sign\si\;
     L([K(\row \xi{\si1}{\si k}),\xi_{\si(k+1)}],
     \xi_{\si(k+2)},\ldots) \\
& + {\tsize\frac{(-1)^{k\ell}}{(k-1)!\,\ell!}}\sum_\si \sign\si\;
     K([L(\row \xi{\si1}{\si \ell}),\xi_{\si(\ell+1)}],
     \xi_{\si(\ell+2)},\ldots) \\
& + {\tsize\frac{(-1)^{k-1}}{(k-1)!\,(\ell-1)!\,2!}}
     \sum_\si \sign \si\; L(K([\xi_{\si 1},\xi_{\si 2}],
     \xi_{\si 3},\ldots), \xi_{\si(k+2)},\ldots) \\
& + {\tsize\frac{(-1)^{(k-1)\ell}}{(k-1)!\,(\ell-1)!\,2!}}
     \sum_\si \sign \si\; K(L([\xi_{\si 1},\xi_{\si 2}],
     \xi_{\si 3},\ldots), \xi_{\si(\ell+2)},\ldots).
\endalign$$
For decomposable tangent bundle valued forms we have the following 
formula for the Fr\"olicher-Nijenhuis bracket in terms of the usual 
operations with vector fields and differential forms, see \cit!{5}, 
or \cit!{9},~8.7. Let $\ph\in\Om^k(M)$, $\ps\in\Om^\ell(M)$, and 
$X,Y\in\X(M)$. Then
$$\align
[\ph\otimes X, \ps\otimes Y] 
     &= \ph \wedge \ps \otimes [X,Y] 
     + \ph \wedge \L_X\ps\otimes Y - 
     \L_Y\ph \wedge \ps \otimes X\tag2 \\
&\qquad + (-1)^k\left(d\ph \wedge i_X\ps \otimes Y 
     + i_Y\ph\wedge  d\ps \otimes X \right).
\endalign$$

\proclaim{\nmb.{4.4}. Proposition} 
Let $M$ be a $\g$-manifold. 
Let $\ph\in\Om_{\text{hor}}^p(M;\g)^\g$ and let 
$\ps\in\Om_{\text{hor}}^q(M;\g)^\g$. Then we have: 
\roster
\item $\ze_\ps$ is horizontal and $\g$-equivariant in the sense that 
       the Lie derivative along any fundamental vector field $\ze_X$ 
       for $X\in\g$ vanishes: 
       $\L_{\ze_X}(\ze_\ps)=[\ze_X,\ze_\ps]=0$. 
\item $[\ze_\ph,\ze_\ps]=-\ze_{[\ph,\ps]^\wedge}$, so 
       $\ze:\Om_{\text{hor}}(M;\g)^\g\to\Om(M;TM)$ is an 
       anti homomorphism of graded Lie algebras.
\endroster
If $\om\in\Om^1(M;\g)^\g$ is a principal connection form with 
principal connection $\Ph=\ze_\om$ and horizontal projection 
$\ch:=Id_{TM}-\Ph$ then we have furthermore:
\roster
\item [3] $[\Ph,\ze_\ps] = -\ze_{d\ps + [\om,\ps]^\wedge}$.
\item $\frac12[\Ph,\Ph] = -\ze_{d\om+\tfrac12[\om,\om]^\wedge}$.
\endroster
\endproclaim

Compare this results with \cit!{9},~11.5, which, however, contains a 
sign mistake in (10). We give here a shorter proof of a stronger 
statement.
  
\demo{Proof} Let $X_1,\dots,X_k$ be a linear basis of the Lie algebra 
$\g$. Then any form $\ph\in\Om^p(M;\g)$ can be uniquely written in 
the form $\ph=\sum_{i=1}^k\ph^i\otimes X_i$ for $\ph^i\in\Om^p(M)$.
Then $\ze_\ph= \sum_{i=1}^k\ph^i\otimes \ze_{X_i}$.

Note that $\ph$ is horizontal if and only if all $\ph^i$ are 
horizontal. Also $\ph$ is $\g$-equivariant, $\ph\in\Om^p(M;\g)^\g$, 
i\.e\. $\L_{\ze_X}\ph=-\ad(X)\o \ph$ for all $X\in\g$, if and only if 
for all $X\in\g$ we have:
$$
\sum_{i=1}^k\L_{\ze_X}\ph^i\otimes X_i = 
     -\sum_{i=1}^k\ph^i\otimes [X,X_i].\tag5
$$
Assertion \therosteritem1 now follows from \thetag5 and 
$$
\L_{\ze_X}\ze_\ph = \sum_{i=1}^k\L_{\ze_X}\ph^i\otimes \ze_{X_i} + 
     \sum_{i=1}^k\ph^i\otimes [\ze_X,\ze_{X_i}].
$$
Using \nmb!{4.3}.(2) we have for general $\ph,\ps\in\Om(M;\g)$
$$\align
[\ze_\ph,\ze_\ps] 
	&= \sum_{i,j}[\ph^i\otimes \ze_{X_i}, \ps^j\otimes \ze_{X_j}] 
= \sum_{i,j}\ph^i \wedge \ps^j \otimes [\ze_{X_i},\ze_{X_j}] \tag6 \\
&\quad + \sum_{i,j}\ph^i \wedge \L_{\ze_{X_i}}\ps^j\otimes \ze_{X_j} 
	- \sum_{i,j}\L_{\ze_{X_j}}\ph^i \wedge \ps^j \otimes \ze_{X_i} \\
&\quad + (-1)^p\sum_{i,j}\left(d\ph^i \wedge i_{\ze_{X_i}}\ps^j 
	\otimes \ze_{X_j} + i_{\ze_{X_j}}\ph^i
	\wedge  d\ps^j \otimes \ze_{X_i} \right).
\endalign$$
If $\ps$ is $\g$-equivariant then from \thetag5 we have
$$
\sum_{i,j} \ph^i \wedge \L_{\ze_{X_i}}\ps^j\otimes \ze_{X_j} 
     =-\sum_{i,j} \ph^i \wedge \ps^j\otimes \ze_{[X_i,X_j]}
     =-\ze_{[\ph,\ps]^\wedge}.
$$ 
So for $\ph$ and $\ps$ both horizontal and $\g$-equivariant 
\thetag6 reduces to assertion \therosteritem2.

If $\ph=\om$, the connection form, then we have
$(-1)^p\sum_{i,j}i_{\ze_{X_j}}\ph^i\wedge  d\ps^j \otimes \ze_{X_i}
     =\sum_{j}d\ps^j \otimes \ze_{X_j} = -\ze_{d\ps}$,
so that \thetag6 reduces to assertion \thetag3.
Similarly, for $\ph=\ps=\om$ formula \thetag6 reduces to 
$[\Ph,\Ph]=-\ze_{[\om,\om]^\wedge + 2d\om}$, so also \thetag4 holds.
\qed\enddemo

\subhead\nmb.{4.5}. Covariant exterior derivative \endsubhead
Let $M$ be a $\g$-manifold of constant rank, let 
$\Ph\in\Om^1(M;TM)^\g$ be a principal  
connection with associated horizontal projection $\ch:= Id_{TM}-\Ph$. 
Let $V$ be any vector space of finite dimension.
Then we define the \idx{\it covariant exterior derivative}
$$
d_\Ph:= \ch^*\o d: \Om^p(M;V)\to \Om^{p+1}_{\text{hor}}(M;V).
$$
We also consider the following mapping as a form of the covariant 
exterior derivative:
$$
\ad(\Ph):=[\Ph,\quad]:\Om^p(M;TM)\to \Om^{p+1}(M;TM).
$$
If a principal connection form $\om:TM\to\g$ exists and 
$\rh:\g\to\frak g\frak l(V)$ is a representation space of $\g$ 
we also consider the following covariant exterior derivative:
$$\gather
d_\om:\Om^p(M;V)\to \Om^{p+1}(M;V)\\
d_\om\Ps := d\Ps + \rh^{\wedge}(\om)\Ps.
\endgather$$

\proclaim{Lemma}
In this situation we have:
\roster
\item $\ad(\Ph)$ restricts to a mapping
     $[\Ph,\quad]:\Om^p_{\text{hor}}(M;\g(M))^\g
     \to \Om^{p+1}_{\text{hor}}(M;\g(M))^\g,$
     where $\g(M)\subset TM$ is the vertical bundle.
\item Let $\rh:\g\to GL(V)$ be a representation. Then
     $d_\Ph$ restricts to a mapping 
     $d_\Ph:\Om^p_{\text{hor}}(M;V)^\g \to 
     \Om^{p+1}_{\text{hor}}(M;V)^\g.$ 
\item For $\ps\in\Om_{\text{hor}}^p(M;\g)^\g$ the two covariant derivatives 
     correspond to each other up to a sign:
     $\ze_{d_\Ph\ps}=-[\Ph,\ze_\ps]$.
\item Let $\rh:\g\to \frak g\frak l(V)$ be a representation. Then
     $d_\om$ restricts to a mapping 
     $d_\om:\Om^p(M;V)^\g \to \Om^{p+1}(M;V)^\g.$ 
     For $\Ps\in \Om_{\text{hor}}^p(M;V)^\g$ and $X\in \g$ we have 
     $i(\ze_X)d_\om\Ps = \rh(\om(\ze_X)-X)\Ps$.
     If $M$ is a free 
     $\g$-manifold, then $d_\om$ also respects horizontality and we 
     have $d_\Ph=d_\om:\Om^p_{\text{hor}}(M;V)^\g \to 
     \Om^{p+1}_{\text{hor}}(M;V)^\g$, where $\Ph=\ze_\om$.
\endroster
\endproclaim

\demo{Proof} \therosteritem1
Let $\Ps\in\Om^p_{\text{hor}}(M;\g(M))^\g$ and $X\in\g$. Then 
formulas \cit!{9},~8.11.(2) give us here
$$\align
&i_{\ze_X}[\Ph,\Ps] = [i_{\ze_X}\Ph,\Ps] 
     + [\Ph,i_{\ze_X}\Ps]  \\
&\qquad -\left(i([\Ph,{\ze_X}])\Ps 
      - (-1)^pi([\Ps,{\ze_X}])\Ph\right) = 0,
\endalign$$
so that $[\Ph,\Ps]$ is again horizontal. It is also $\g$-equivariant 
since $[\ze_X,[\Ph,\Ps]]=[[\ze_X,\Ph],\Ps]+[\Ph,[\ze_X,\Ps]]=0$ for all 
$X\in\g$ by the graded Jacobi identity. That it has vertical values 
can be seen by contemplating one of the formulas in \nmb!{4.3}.

\therosteritem2 Let $\Ps\in\Om^p_{\text{hor}}(M;V)^\g$. For 
$X\in\g$ we have $\L_{\ze_X}\Ps=-\rh(X)\o \Ps$, then $d_\Ph\Ps$ is 
again $\g$-equivariant, since we have
$$\align
\L_{\ze_X}\ch^*d\Ps &= \L_{\ze_X}(d\Ps\o \La^{p+1}\ch) 
     = \L_{\ze_X}(d\Ps)\o \La^{p+1}\ch + 
     d\Ps\o \L_{\ze_X}(\La^{p+1}\ch)  \\
&= (d\L_{\ze_X}\Ps)\o \La^{p+1}\ch + 0 = \ch^*d(-\rh(X)\o\Ps)
     = -\rh(X)\o(\ch^*d\Ps),
\endalign$$
and clearly horizontal.

\therosteritem3 
Let again $X_1,\dots,X_k$ be a linear basis of the Lie algebra 
$\g$ and consider 
$\ps=\sum_{i=1}^k\ps^i\otimes X_i \in \Om_{\text{hor}}^p(M;\g)^\g$ 
for $\ps^i\in\Om^p_{\text{hor}}(M)$. Then we use \cit!{9},~8.7.(5) to 
get 
$$\align
(-1)^{p+1}[\Ph,\ze_\ps] &= [\ze_\ps,\Ph] 
     = \sum_i [\ps^i\otimes \ze_{X_i},\Ph] \\
&= \sum_i \left(\ps^i\wedge [\ze_{X_i},\Ph] 
     - (-1)^p \L_{\Ph}\ps^i\otimes\ze_{X_i} +(-1)^p 
     d\ps^i\wedge i(\ze_{X_i})\Ph \right).
\endalign$$
Since $\Ph$ is $\g$-equivariant we have $[\ze_{X_i},\Ph]=0$. Moreover 
we have 
$\L_\Ph\ps^i = i_\Ph d\ps^i - d i_\Ph\ps^i = i_\Ph d\ps^i - 0$ and 
$i(\ze_{X_i})\Ph = \ze_{X_i}$.
Thus we get
$$\align
(-1)^{p+1}[\Ph,\ze_\ps] &= (-1)^p \sum_i 
     (d\ps^i - i_\Ph d\ps^i)\otimes\ze_{X_i}\\
&= (-1)^p \sum_i \ch^*d\ps^i\otimes\ze_{X_i} = (-1)^p 
     \ze_{d_\Ph\ps}. 
\endalign$$

\therosteritem4
Let $\Ps\in\Om^p(M;V)^\g$. For 
$X\in\g$ we have $\L_{\ze_X}\Ps=-\rh(X)\o \Ps$, then $d_\om\Ps$ is 
again $\g$-equivariant, since we have
$$\align
\L_{\ze_X}(d\Ps + \rh^{\wedge}(\om)\Ps) &= d\L_{\ze_X}\Ps + 
     \rh^{\wedge}(\L_{\ze_X}\om)\Ps + \rh^{\wedge}(\om)\L_{\ze_X}\Ps\\
&= -d\rh(X)\Ps  
    - \rh^{\wedge}([X,\om]^\wedge)\Ps - \rh^{\wedge}(\om)\rh(X)\Ps\\
&= -\rh(X)(d\Ps + \rh^{\wedge}(\om)\Ps).
\endalign$$
For $\Ps\in \Om_{\text{hor}}^p(M;V)^\g$ and $X\in \g$ we use
$i_{\ze_X}\Ps=0$ and $\L_{\ze_X}\Ps= -\ad(X)\Ps$ to get
$$\align
i_{\ze_X}(d\Ps + \rh^\wedge(\om)\Ps) &= i_{\ze_X}d\Ps + 
	di_{\ze_X}\Ps + \rh(i_{\ze_X}\om)\Ps - \rh^\wedge(\om)i_{\ze_X}\Ps\\
&= \L_{\ze_X}\Ps + \rh(\om(\ze_X))\Ps =  \rh(\om(\ze_X)-X)\Ps.
\endalign$$
Let now $M$ be a free $\g$-manifold then $\om(\ze_X)-X=0$ and 
$d_\om\Ps$ is again horizontal. 
We use the principal connection $\Ph$ 
to split each vector field into the sum of a horizontal one 
and a vertical one. 
If we insert one vertical vector field, say $\ze_X$ for 
$X\in \frak g$, into $d_\Ph\Ps-d_\om\Ps$, we get 0.
Let now all vector fields $\xi_i$ be horizontal, then we get
$$\gather
(d_\Ph\Ps)(\row \xi 0k) = (\ch^*d\Ps)(\row \xi 0k) = 
	d\Ps(\row \xi0k),\\
(d\Ps + [\om,\Ps]^\wedge)(\row\xi0k) = d\Ps(\row \xi0k).\qed
\endgather$$
\enddemo

\subhead\nmb.{4.6}. Curvature \endsubhead
Let $M$ be a $\g$-manifold.
If there exists a principal connection  $\Ph$ then this is a 
projection onto the integrable vertical distribution induced by $\g$, 
and the formula \nmb!{4.3}.\thetag1 
for the Fr\"olicher-Nijenhuis bracket reduces to
$$
R(\xi,\et)=\frac12[\Ph,\Ph](\xi,\et)=\Ph[\xi-\Ph\xi,\et-\Ph\et].
$$
$R\in\Om^2_{\text{hor}}(M;TM)^\g$ is called the \idx{\it curvature} of 
the connection $\Ph$. From the graded Jacobi identity of the 
Fr\"olicher-Nijenhuis bracket we get immediately the \idx{\it Bianchi 
identity}
$$
[\Ph,R]=\tfrac12[\Ph,[\Ph,\Ph]]=0.
$$
Note that the kernel of $\ad(R)$ is invariant under $\ad(\Ph)$, and 
$\ad(\Ph)^2=0$ on it. It gives rise to a cohomology, depending on 
$\Ph$.

If $\om\in\Om^1(M,\g)^\g$ is a principal connection form, then 
formula \thetag4 in proposition \nmb!{4.4} suggests to define
$$
\Om:= d\om+\tfrac12[\om,\om]^\wedge
$$
as the \idx{\it curvature form} of $\om$; so we have $R=-\ze_\Om$.
Then \nmb!{4.4},~\thetag3 
suggests that the \idx{\it Bianchi identity} should have 
the form 
$d_\om\Om= d\Om+[\om,\Om]^\wedge = 0$.
Indeed this is true and it follows directly 
from the graded Jacobi identity in 
$(\Om(M,\g),[\quad,\quad]^\wedge)$. 

\proclaim{\nmb.{4.7}. Proposition}
Let $M$ be a $\g$-manifold with 
principal connection $\Ph$ and horizontal projection 
$\ch:=Id_{TM}-\Ph$. Then we have:
\roster
\item $d_\Ph\o \chi^*-d_\Ph =\ch^*[d,\ch^*]= 
       \chi^*\o i_R:\Om^p(M;V)\to\Om_{\text{hor}}^{p+1}(M;V)$, 
       where $R$ is the curvature and $V$ is any vector space and 
       $i_R$ is the insertion operator.
\item $d_\Ph\o d_\Ph = \chi^*\o i(R)\o d:
       \Om^p(M;V)\to\Om_{\text{hor}}^{p+1}(M;V)$. 
\item If $\om\in\Om^1(M;\g)^\g$ is a principal connection form then 
       the curvature form $\Om=d\om+\tfrac12[\om,\om]^\wedge$ 
       satisfies
       $i(\ze_X)\Om = [\om(\ze_X)-X,\om]^\wedge - d(\om(\ze_X))$. If 
       $M$ is a free $\g$-manifold then $\Om$ is horizontal and 
       $\Om = d_\Ph\om = d_\om\om\in\Om^2_{\text{hor}}(M;\g)^\g$.
\endroster
\endproclaim

Note that by \therosteritem2 the kernel of $\chi^*\o i(R)\o d$ is 
invariant under $d_\Ph$, which gives rise to a cohomology associated 
to it.

\demo{Proof}
\therosteritem1 
For $\Ps\in\Om(P;V)$ we have
$$\allowdisplaybreaks\align
(d_\Ph\chi^*\Ps)&(\row \xi0k) 
     = (d\chi^*\Ps)(\chi(\xi_0),\dots,\chi(\xi_k))\\
&=\sum_{0\leq i\leq k}(-1)^i\chi(\xi_i)((\chi^*\Ps)
     (\chi(\xi_0),\dots,\widehat{\chi(\xi_i)},\dots,\chi(\xi_k)))\\
&\quad+\sum_{i<j}(-1)^{i+j}(\chi^*\Ps)
     ([\chi(\xi_i),\chi(\xi_j)],\chi(\xi_0),\dotsc\\
&\qquad\qquad\qquad\qquad\dots,\widehat{\chi(\xi_i)},\dots,
     \widehat{\chi(\xi_j)},\dotsc)\\
&=\sum_{0\leq i\leq k}(-1)^i\chi(\xi_i)(\Ps
     (\chi(\xi_0),\dots,\widehat{\chi(\xi_i)},\dots,\chi(\xi_k)))\\
&\quad+\sum_{i<j}(-1)^{i+j}\Ps([\chi(\xi_i),\chi(\xi_j)]
     -\Ph[\chi(\xi_i),\chi(\xi_j)],\chi(\xi_0),\dotsc\\
&\qquad\qquad\qquad\qquad\dots,\widehat{\chi(\xi_i)},\dots,
     \widehat{\chi(\xi_j)},\dotsc)\\
&=(d\Ps)(\chi(\xi_0),\dots,\chi(\xi_k)) +
     (i_{R}\Ps)(\chi(\xi_0),\dots,\chi(\xi_k))\\
&=(d_\Ph + \chi^*i_{R})(\Ps)(\row \xi0k).
\endalign$$

\therosteritem2
$d_\Ph d_\Ph = \chi^*d\chi^*d = (\chi^*i_R + \chi^*d)d = \chi^*i_R d$ 
holds by \therosteritem1.

\therosteritem3
For $X \in \frak g$ we have 
$$\multline
i_{\ze_X}(d\om + \frac12[\om,\om]^\wedge ) = i_{\ze_X}d\om  +
     \frac 12 [i_{\ze_X}\om,\om]^\wedge  
     - \frac 12 [\om,i_{\ze_X}\om]^\wedge  = \\
= \Cal L_{\ze_X}\om -d(\om(\ze_X)) + [\om(\ze_X),\om]
=  [\om(\ze_X)-X,\om] - d\om(\ze_X).
\endmultline $$
If $M$ is a free $\g$-manifold then this is zero, and on horizontal 
vectors $d_\Ph\om$ and $d_\om\om$ coincide.
\qed\enddemo

\head\totoc\nmb0{5}. Homogeneous $\g$-manifolds \endhead

\subhead\nmb.{5.1}. Homogeneous free $\g$-manifolds and Maurer-Cartan 
forms \endsubhead
Recall that a $\g$-valued 1-form $\ka$ on a manifold $M$ is called 
\idx{\it Maurer-Cartan form} if $\ka_x:T_xM\to \g$ is a linear 
isomorphism for each $x\in M$ and if $\ka$ satisfies the 
\idx{\it Maurer-Cartan equation} $d\ka + \tfrac12[\ka,\ka] =0$.
This concept is also sometimes called a
\idx{\it flat Cartan connection}, and a manifold with a flat Cartan 
connection is sometimes called a \idx{\it principal homogeneous 
space}. See \cit!{7} for Maurer-Cartan forms. 

\proclaim{Lemma}
To each free transitive $\g$-action $\ze:\g\to \X(M)$ 
there corresponds a unique Maurer-Cartan form $\ka:TM\to\g$, given by 
$\ka_x=\ze_x\i$, and conversely. Then $\ka$ is $\g$-equivariant with 
respect to the $\g$-action $\ze:\g\to \X(M)$, and $\ka$ is the unique 
principal connection form on the $g$-manifold $M$.
\endproclaim

Note also that an action of a Lie algebra $\g$ is free if and only if 
the associated pseudogroup has discrete isotropy groups.

\demo{Proof}
If $\ka\in \Om^1(M;\g)$ and $\ze:\g\to \X(M)$ are inverse to each 
other then for $X,Y\in\g$ we have 
$$\align
(d\ka+\tfrac12[\ka,\ka]^\wedge)(\ze_X,\ze_Y) &= \ze_X(\ka(\ze_Y)) - 
     \ze_Y(\ka(\ze_X)) - \ka([\ze_X,\ze_Y]) + [\ka(\ze_X),\ka(\ze_Y)]\\
&= - \ka([\ze_X,\ze_Y]) + [X,Y] \\
&= -\ka([\ze_X,\ze_Y]-\ze_{[X,Y]}),
\endalign$$
so that $\ze$ is a Lie algebra homomorphism if and only if $\ka$ 
fulfills the Maurer-Cartan equation.
For fixed $\g$-action $\ze$ the form $\ka$ is $\g$-equivariant, since 
we have
$\L_{\ze_X}\ka = i_{\ze_X}d\ka + di_{\ze_X}\ka 
= - i_{\ze_X}(\tfrac12[\ka,\ka]^\wedge) + dX
= - [i_{\ze_X}\ka,\ka] + 0 = -\ad(X)\ka$.
Thus $\ka$ is a principal connection form for this $\g$-action, and 
it is the unique one by proposition \nmb!{5.7} below.
\qed\enddemo

\subhead\nmb.{5.2}. Cartan's developing \endsubhead
It is well known that a free homogeneous $G$-manifold may be 
identified with the Lie group $G$ by fixing a point. For 
$\g$-manifolds the situation is more complicate. The following result 
may also be found in \cit!{7}.

\proclaim{Proposition}
Let $M$ be a free transitive $\g$-manifold which is simply connected. 
Let $G$ be a Lie group with Lie algebra $\g$. Then there exist 
$\g$-equivariant local diffeomorphisms $M\to G$. Namely for each 
$x\in M$ there is a unique $\g$-equivariant smooth mapping 
$C_x:M\to W_x$ with $C_x(x)=e$ which is locally a diffeomorphism, 
where $W_x\subset G$ is defined in \nmb!{2.3}.

If the $\g$-action on $M$ integrates to a $G$-action on $M$, then 
this mapping is automatically a diffeomorphism.
\endproclaim

The embedding $M\to G$
is called \idx{\it Cartan's developing}. Its origins 
lie in Cartan's developing of a locally Euclidean space into the 
standard Euclidean space. If Cartan's developing is injective then 
the $\g$-manifold $M$ admits an extension to a $\g$-manifold which is 
isomorphic to the Lie group $G$ with the left action of $\g$.

\demo{Proof}
We consider again $G\x M$ with the graph foliation as in \nmb!{2.3}. 
Then $\operatorname{pr}_2:G\x M \to M$ is a principal $G$-bundle 
with left multiplication as principal action, and since $M$ is a 
free $\g$-manifold the graph foliation is transversal to the fibers 
of $\operatorname{pr}_2$ and is the horizontal foliation of a 
principal connection on the $G$-bundle. For $x\in M$ the restriction 
of $\operatorname{pr}_2$ to the leaf $L(x)$ through $(e,x)$ is a 
$\g$-equivariant covering mapping which is a diffeomorphism since $M$ 
is simply connected. Then
$$
\ph: M @>{\operatorname{pr}_2\i}>> L(x) 
     @>{\operatorname{pr}_1}>> W_{x}\subset G
$$
is the looked for $\g$-equivariant mapping which locally is a 
diffeomorphism since $\operatorname{pr}_1$ also is one. 
                                                         
It remains to show that $\ph$ is a diffeomorphism if the $\g$-action 
is complete. We consider the $\g$-action on $G\x M$ on the factor $M$ 
alone in this case. Then the graph foliation gives us a flat principal 
connection for this action, see section \nmb!{3}, 
and by proposition \nmb!{6.6} below this 
connection is complete. Thus $\operatorname{pr}_1:L(x)\to G$ is also 
a covering map, and since $G$ is simply connected it is a 
diffeomorphism also and we are done.
\qed\enddemo

\subhead\nmb.{5.3}. Example \endsubhead 
The result of proposition \nmb!{5.2} is the best possible in general, 
as the following example shows. 
Let $G$ be a simply connected Lie group, let $W$ be a not simply 
connected open subset of $G$, and let $M$ be a simply connected 
subset of the universal cover of $W$ such that the projection 
$\ph:M\to W$ is still surjective. We have an action of the Lie 
algebra $\g$ of $G$ on $M$ by pulling back all left invariant vector 
fields on $G$ to $M$ via $\ph$. Then $\ph$ is as constructed in 
\nmb!{5.2}, but it is only locally a diffeomorphism.

For example, let $W$ be an annulus in $\Bbb R^2$, and let $M$ be a 
piece of finite length of the spiral covering the annulus.
Other examples can be found in \cit!{10}, \cit!{11}.

\subhead\nmb.{5.4} \endsubhead
As an immediate application of the Cartan developing, we have the 
following proposition:

\proclaim{Proposition}
Let $H$ be a connected linear Lie group of finite type, let $G$ be 
the simply connected full prolongation of $H$ such that $G/H$ is the 
standard maximally homogeneous $H$-structure (see \cit!{1}).

Then for any manifold $N$ with a locally flat $H$-structure 
$p:M\to N$ there exists a map $\ps:N\to G/H$, which is a local 
isomorphism of $H$-structures.
\endproclaim

In the case of a flat conformal structure we obtain the well known 
developing of a locally flat conformal manifold into the conformal 
sphere.

\demo{Proof}
The mapping $\ps$ is 
the unique one making the following diagram commutative:
$$\CD
M^\infty @>\ph>> G \\
@V{p_\infty}VV @VV{p_0}V \\
N @>\ps>> G/H. 
\endCD$$
Here $M^\infty$ is the full prolongation of the $H$-structure $p$ 
with the natural free transitive action of the Lie algebra $\g$ of 
$G$, see \cit!{1}, and $\ph$ is the Cartan developing of the 
$\g$-manifold $M^\infty$ into $G$.
\qed\enddemo

\subhead\nmb.{5.5}. The dual $\g$-action for simply connected 
homogeneous free $\g$-ma\-ni\-folds \endsubhead 

As motivation we recall that on a Lie group $G$ (viewed as a 
homogeneous free right $G$-manifold) the fundamental vector fields 
correspond to the left invariant ones; they generate right 
translations, and correspond to the left Maurer-Cartan form $\ka$ on 
$G$. The diffeomorphisms which commute with all right translations 
are exactly the left translations; the vector fields commuting 
with all left invariant ones are exactly the right invariant ones; 
they generate left translations, and correspond to the right 
Maurer-Cartan form $\hat\ka= \Ad.\ka$. 

Now let $M$ be a free homogeneous $\g$-manifold with action 
$\ze:\g\to \X(M)$ and the corresponding principal connection form 
$\ka=\ze\i$, see \nmb!{5.1}. Let $G$ be a Lie group with Lie algebra 
$\g$. Choose a point $x_0\in M$. 
Assume that $M$ is simply connected and consider the Cartan developing 
$C_{x_0}:M\to G$. Then for $X\in \g$ the fundamental vector field 
$\ze_X\in \X(M)$ is $C_{x_0}$-related to the left invariant vector 
field $L_X$ on $G$. Let now $\hat\ze_X\in\X(M)$ denote the unique 
vector field on $M$ which is $C_{x_0}$-related to the right invariant 
vector field $R_X\in\X(G)$. Since $[L_X,R_Y]=0$ we get 
$[\ze_X,\hat\ze_Y]=0$, and even each local vector field $\xi\in\X(U)$ 
for connected open $U\subset M$
with $[\xi,\ze_Y]=0$ for all $Y\in\g$ extends to one of the form 
$\hat\ze_X$. So we get a Lie algebra anti homomorphism
$$
\hat\ze:\g\to \X(M)
$$
whose image is the centralizer algebra
$$
Z_{\X(M)}(\g):=\{\et\in\X(M):[\et,\ze_X]=0 \text{ for all }X\in\g\}
$$
This (`right') action $\hat\ze$ of $\g$ on $M$ which commutes with 
the original action $\ze$ is called the \idx{\it dual action}. We 
have also the dual principal connection form $\hat\ka$, inverse to 
$\hat\ze$, see \nmb!{5.1}.

Note that for a free homogeneous $\g$-manifold $M$ which is not simply 
connected, the dual action of $\g$ does not exist in general and the 
centralizer algebra $Z_{\X(M)}(\g)$ is smaller than $\g$: As an 
example we consider a Lie group $G$ with the right action on 
$H\backslash G$ for a discrete subgroup $H$ which is not central. 
Then the associated action of the Lie algebra $\g$ is free, but its 
centralizer $Z_{\X(H\backslash G)}(\g)$ is isomorphic to the 
subalgebra $\g^H:=\{X\in \g:\Ad(h)X=X\text{ for all }h\in H\}$.

\subhead\nmb.{5.6}. Homogeneous $G$-manifolds \endsubhead
As a motivation for what follows we consider here homogeneous 
$G$-manifolds. So let $G$ be a connected Lie group with 
Lie algebra $\g$, multiplication $\mu:G\x G\to G$, and for $g\in G$ 
let $\mu_g, \mu^g:G\to G$ denote the left and right translation, 
$\mu(g,h)=g.h=\mu_g(h)=\mu^h(g)$. 
Let $H\subset G$ be a closed subgroup with Lie algebra $\h$.

We consider the right coset space $M=H\backslash G$, the canonical 
projection $p:G\to H\backslash G$, the initial point 
$o=p(e)\in H\backslash G$ and the canonical right action of $G$ on 
the right coset space $H\backslash G$, denoted by 
$\bar\mu^g:H\backslash G \to H\backslash G$. Then for $X\in\g$ 
the left invariant vector field $L_X\in\X(G)$ is $p$-related to the 
fundamental vector field $\ze_X\in\X(H\backslash G)$ of $\bar\mu^g$.

Suppose now that we are given a principal connection form 
$\om\in \Om^1(H\backslash G;\g)^\g$. Then $\L_{\ze_X}\om=-\ad(X)\om$ 
implies in turn
$$\gather
(\bar\mu^{\exp(tX)})^*\om = (\Fl^{\ze_X}_t)^*\om = 
     e^{-\ad(tX)}\om = \Ad(\exp(-tX))\om\\
(\bar\mu^g)^*\om = \Ad(g\i)\om \quad\text{ for all }g\in G\tag1\\
\om_{o.g} = \Ad(g\i)\o \om_o\o T_o(\bar\mu^{g\i})\tag2
\endgather$$
We also get a reductive decomposition of the Lie algebra $\g$ as
$$
\g=\h \oplus \om_o(T_oH\backslash G) =: \h \oplus \frak m_o,
\tag3$$
where $\frak m_o$ is a linear complement to $\h$ which is invariant 
under $\Ad(H)$.

Conversely any $\Ad(H)$-invariant linear complement $\frak m_o$ of 
$\h$ in $\g$ defines a principal connection form on the $\g$-manifold 
$H\backslash G$ as follows: we consider the $H$-equivariant linear 
mapping
$$
\om_o: T_o(H\backslash G) @<{\ze_o}<{\cong}< 
     \g/\h\cong \frak m_o \subset \g
$$
and extend it to a principal connection form $\om$ by \thetag2.

There is a bijective correspondence between principal connection 
forms $\om\in\Om^1(H\backslash G,\g)^\g$ and principal connections 
$\tilde\om\in\Om^1(G;\h)^H$ on the principal fiber bundle 
$p:G\to H\backslash G$ with {\it left} principal action of $H$, 
which is given by 
$$\align
\tilde\om_g :&= T_g(\mu^{g\i}) - 
     \Ad(g)\o\om_{p(g)}\o T_gp:T_gG\to \h, \tag 4\\
\tilde\om :&= \ka - \Ad. p^*\om, \tag5
\endalign$$
where $\ka$ denotes the right Maurer-Cartan form.
It is easily checked that $\tilde\om$ is a principal connection for 
$p:G\to H\backslash G$: since the principal $H$-action is the left 
action on $G$ we have $(\mu^h)^*\tilde\om = \Ad(h)\om$, and $\om$ 
reproduces the generators in $\h$ of right invariant vector fields on 
$G$. The principal curvature of $\tilde\om$ is given by 
$d\tilde\om-\tfrac12[\tilde\om,\tilde\om]_\h^\wedge$: see \cit!{9}, 
proof of 11.2.(3); compare with \cit!{8},~I,~Chap\.~X.
The curvature form of $\tilde\om$ is related to the curvature form 
$\Om=d\om+\tfrac12[\om,\om]_\g^\wedge$ of $\om$ by 
$$\align
d\tilde\om-\tfrac12[\tilde\om,\tilde\om]_\h^\wedge 
&= d\ka - \tfrac12[\ka,\ka]^\wedge_\g - d\Ad \wedge p^*\om - 
     \Ad. p^*d\om \\
&\quad + [\ka,\Ad. p^*\om]^\wedge_\g - 
     \tfrac12[\Ad. p^*\om,\Ad. p^*\om]^\wedge_\g\\
&= - \Ad. p^*(d\om + \tfrac12[\om,\om]^\wedge_\g) 
     = - \Ad. p^*\Om,
\endalign$$ 
since for the right Maurer-Cartan form $\ka$ the Maurer-Cartan 
equation is given by $d\ka-\tfrac12[\ka,\ka]_\g^\wedge=0$, 
and since for $X\in\g$ we have:
$$\align
d\Ad (T(\mu^g)X) &= \ddt \Ad(\exp(tX).g) = \ad(X)\Ad(g) 
     = \ad(\ka(T(\mu^g)X))\Ad(g),\\
d\Ad &=(\ad\o\ka)\Ad. 
\endalign$$

\proclaim{\nmb.{5.7} Proposition}
Let $M$ be a homogeneous $\g$-space. Then there exists a unique 
principal connection $\Ph=Id$ on $M$. 

On the other hand let $M$ be an effective homogeneous $\g$-space. 
Then principal 
connection forms $\om$ correspond to reductive decompositions 
$\g=\g_x + \frak m_x$, where $\g_x$ is the isotropy subalgebra of a 
point $x\in M$, and where $\frak m_x$ is an $\Ga(\g)_x$-invariant 
complementary subspace. 
\endproclaim

\demo{Proof} The first statement is obvious.
 
We first check that for an effective homogeneous $\g$-manifold $M$ 
the homomorphism 
$\operatorname{germ}_x\o\ze^M:\g\to \X(M)_{\text{germs at }x}$ is 
injective for each $x\in M$. Let $\frak k_x$ denote its kernel. Since 
$\g$ is finite dimensional, we have $\frak k_y=\frak k_x$ for $y$ 
near $x$, and since $M$ is connected, this holds even for all 
$y\in M$. So $\frak k_x$ is in the kernel of $\ze^M:\g\to \X(M)$ 
which is zero since the $\g$-action on $M$ is supposed to be 
effective.

As in \nmb!{2.3} we consider the graph foliation of the 
$\g$-manifold $M$ on $G\x M$, where $G$ is a connected Lie group with 
Lie algebra $\g$, the leaf $L(x)$ through $(e,x)$ of it, and the 
open set $W_{x}=\operatorname{pr}_1(L(x))\subset G$. 

Consider a principal connection form $\om\in\Om^1(M;\g)^\g$.
Then $\L_{\ze_X}\om=-\ad(X)\om$ 
implies 
$(\Fl^{\ze_X}_t)^*\om = 
     e^{-\ad(tX)}\om = \Ad(\exp(-tX))\om$, this holds then for all 
elements of the pseudogroup $\Ga(\g)$ of the form \nmb!{2.2}.(5) and 
finally we get for each smooth curve $c:[0,1]\to W_x$ with $c(0)=e$ 
which is liftable to $L(x)$:
$$
(\ga_x(c))^*\om = \Ad(c(1)\i)\om.
\tag1$$
Thus we get a reductive decomposition of  the Lie algebra $\g$ as
$$
\g=\g_x \oplus \om_x(T_xM) =: \g_x \oplus \frak m_x,\tag2
$$
where $\frak m_x$ is a linear complement to $\g_x$ which is invariant 
under $\Ad(\Ga(\g)_x)$, see also \nmb!{2.6}. 

If conversely we are given a reductive decomposition as in \thetag2 
which is invariant under $\Ad(\Ga(\g)_x)$, then we consider the 
$\Ga(\g)_x$-equivariant linear mapping
$\om_x:T_xM @<{\ze_x}<\cong< \g/\g_x \cong \frak m_x \subset \g$
and we use \thetag1 to define 
$\om\in\Om^1(M;\g)^\g$ by
$$
\om_{\ga_x(c)(x)} = \Ad(c(1)\i)\o\om\o T_x(\ga_x(c)\i),
\tag3$$
for each smooth curve $c:[0,1]\to W_x$ with $c(0)=e$ which is 
liftable to $L(x)$. Since any element of the pseudogroup $\Ga(\g)$ is 
of this form (see \nmb!{2.3}) we get a well defined principal 
connection form on $M$.
\qed\enddemo

\proclaim{\nmb.{5.8}. Theorem}
Let $M$ be a homogeneous effective $\g$-space. Let $x_0\in M$ and let 
$P_{x_0}(\g)$ be the set of all germs at $x_0$ of transformations in 
$\Ga(\g)$. 

Then $\ev_{x_0}:P_{x_0}(\g) \to M$ is the projection of a 
smooth principal fiber bundle with structure group the isotropy group 
$\Ga(\g)_{x_0}$ (see \nmb!{2.6}) and with principal right action just 
composition from the right, and its smooth structure is the unique 
one for which the smooth curves $[0,1]\to P_{x_0}(\g)$ correspond 
exactly to the germs along $[0,1]\x \{x_0\}$ of smooth mappings 
$\ph:[0,1]\x U\to M$ with $\ph_t\in\Ga(\g)$ for all $t$.
The total space $P_{x_0}(\g)$ is connected. 

We have a canonical free transitive $\g$-action 
$\ze^P:\g\to\X(P_{x_0}(\g))$ which is given by
$$
\ze^P_X(\ph):= \ze^M_X\o\ph = T\ph\o \ze^M_{\Ad(\ph\i)X},
$$
and whose corresponding pseudogroup is generated by the local flows 
$\ph\mapsto \Fl^{\ze^M_X}_t\o \ph$. Each vector field $\ze^P_X$ is 
invariant under the pullback  by the principal right action of 
$\Ga(\g)_{x_0}$.
The projection $\ev_{x_0}:P_{x_0}(\g)\to M$ is 
$\g$-equivariant, so the vector fields $\ze^P_X$ and $\ze^M_X$ are 
$\ev_{x_0}$-related.
Its associated Maurer-Cartan form is called $\ka$.
                                           
There exists also the dual free $\g$-action (a Lie algebra anti 
homomorphism) $\hat\ze^P:\g\to \X(P_{x_0}(\g))$, given by
$$
\hat\ze^P_X:= T\ph\o \ze^M_X = \ze^M_{\Ad(\ph)X}\o \ph = 
     \ze^P_{\Ad(\ph)X}(\ph).
$$
Its associated Maurer-Cartan form is called $\hat\ka$, and
its corresponding pseudogroup is generated by 
$\ph\mapsto \ph\o \Fl^{\ze^M_X}_t$.  
The pullback of $\hat\ze^P_X$ by the principal right action of 
$\ps\in\Ga(\g)_{x_0}$ is given by 
$(\ps^*)^*\hat\ze^P_X = \hat\ze^P_{\Ad(\ps\i)X}$.

The principal connections forms $\om\in\Om^1(M;\g)^\g$ 
correspond bijectively to principal connections forms 
$\tilde\om\in\Om^1(P_{x_0}(\g),\g_{x_0})^{\Ga(\g)_{x_0}}$ on the 
principal $\Ga(\g)_{x_0}$-bundle $P_{x_0}(\g)\to M$ via 
$$
\tilde\om_\ph = - \hat\ka_\ph - \Ad(\ph\i)\o\om_{\ph(x_0)}\o 
     T_\ph(\ev_{x_0}):T_\ph(P_{x_0}(\g))\to \g_{x_0}.
$$
The principal curvature forms are then related by 
$$
\tilde\Om:= d\tilde\om+\tfrac12[\tilde\om,\tilde\om]^\wedge_{\g_{x_0}} = 
     -(\Ad\o \operatorname{inv}).
     \ev_{x_0}^*(d\om + \tfrac12[\om,\om]^\wedge_\g) = 
     -(\Ad\o \operatorname{inv}).\ev_{x_0}^*\Om.
$$
\endproclaim

The manifold $P_{x_0}(M)$ is not simply connected in general (e\.g\. 
a Lie group); nevertheless the dual action is defined.

\demo{Proof}
Recall first from the proof of \nmb!{5.7} that for an effective 
homogeneous $\g$-manifold $M$ the homomorphism 
$\operatorname{germ}_x\o\ze^M:\g\to \X(M)_{\text{germs at }x}$ is 
injective for each $x\in M$. 

As in \nmb!{2.3} we consider the graph foliation of the 
$\g$-manifold $M$ on $G\x M$, where $G$ is a connected Lie group with 
Lie algebra $\g$, the leaf $L(x_0)$ through $(e,x_0)$ of it, and the 
open set $W_{x_0}=\operatorname{pr}_1(L(x_0))\subset G$. 

Now we choose a splitting $\g=\g_{x_0}\oplus\frak m$, where $\frak m$ 
is a linear complement to the isotropy algebra $\g_{x_0}$. 
Let us now consider a small open ball $B$ with center 0 in $\g$, 
its diffeomorphic image $\exp(B)=:W\subset W_{x_0}\subset G$, 
and $\tilde W\subset L(x_0)$, an 
open neighborhood of $(e,x_0)$ in $L(x_0)$ such that 
$\operatorname{pr}_1|\tilde W: \tilde W \to W$ is a 
diffeomorphism. Let $B^{\frak m}:= B\cap \frak m$, 
$W^{\frak m}:=\exp(B^{\frak m})$, and $\tilde W^{\frak m}:=
(\operatorname{pr}_1|\tilde W)\i(W^{\frak m})$, and choose now $B$ so 
small that $\operatorname{pr_2}:\tilde W^{\frak m}\to M$ is a 
diffeomorphism onto an open neighborhood $U$ of $x_0$. We consider 
the composed diffeomorphism 
$$
\la:U @>{\operatorname{pr}_2\i}>> \tilde W^{\frak m} 
     @>{\operatorname{pr}_1}>> W^{\frak m} @>{\exp\i}>> B^{\frak m}.
$$ 
Now for $f\in P_{x_0}(\g)$ with 
$\ev_{x_0}(f)=f(x_0)\in U$ we define 
$$
\ph(f) = \ga_{x_0}(c_{f(x_0)})\i\o f \in \Ga(\g)_{x_0},
$$
where $c_{f(x_0)}:[0,1]\to W^{\frak m}$ is the curve 
$c_{f(x_0)}(t)=\exp(t\la({f(x_0)}))$.  

Next we choose an open cover $(U_\al)$ of $M$ with transformations 
$f_\al: V_\al \to U_\al$ in the pseudogroup $\Ga(\g)$, where $V_\al$ 
is a connected open neighborhood of $x_0$ in $U$, and we define
$$\gather
\ph_\al:\ev_{x_0}\i(U_\al)= 
     P_{x_0}(\g)|U_\al \to U_\al \x \Ga(\g,x_0)\\
\ph_\al(f):= (f(x_0),\ph(f_\al\i\o f) )=(f(x_0), 
     \ga_{x_0}(c_{f_\al\i f(x_0)})\i\o f_\al\i\o f).
\endgather$$
These give a smooth principal fiber bundle atlas for $P_{x_0}(\g)$ 
since for $(x,h)\in (U_\al\cap U_\be) \x \Ga(\g)_{x_0}$ we have 
$$
\ph_\al\ph_\be\i(x,h) = (x,\ga_{x_0}(c_{f_\al\i(x)})\i\o f_\al\i\o 
     f_\be\o \ga_{x_0}(c_{f_\be\i(x)})\o h).
$$ 
The smooth structure on $P_{x_0}(\g)$ induced by this atlas 
is the unique one where the smooth curves are exactly as described in 
the theorem, since this is visibly the case in each chart.
Thus by the lemma in \nmb!{2.3} the total space $P_{x_0}(\g)$ is 
connected.

For $X\in\g$ and $\ph\in P_{x_0}(\g)$ we have 
$\tfrac \partial{\partial t}\Fl^{\ze^M_X}_t\o \ph = 
\ze^M_X\o\Fl^{\ze^M_X}_t\o\ph$, 
so a smooth vector field on $P_{x_0}(\g)$ is defined by 
$$
\ze^P_X(\ph):= \ze^M_X\o\ph = T\ph\o T\ph\i\o\ze^M_X\o\ph = 
     T\ph\o (\ph^*\ze^M_X) = T\ph\o \ze^M_{\Ad(\ph\i)X},
$$
where we used \nmb!{2.4}, 
and its local flow is given by
$\Fl^{\ze^P_X}_t(\ph) = \Fl^{\ze^M_X}_t\o\ph$. Clearly 
$\ze^P:\g\to \X(P_{x_0}(\g))$ is a $\g$-action, which is free, since 
for each $x\in M$ the homomorphism
$\operatorname{germ}_x\o\ze^M:\g\to \X(M)_{\text{germs at }x}$ is 
injective. 
Consider now $\ps\in\Ga(\g)^{x_0}$ and its principal right action 
$\ps^*$ on $P_{x_0}(\g)$; it acts trivially by pullback on each 
vector field $\ze^P_X$ since we have:
$$\align
((\ps^*)^*\ze^P_X)(\ph) &= (T(\ps^*)\i\o\ze^P_X\o \ps^*)(\ph) 
     = (\ps\i)^*(\ze^P_X(\ph\o\ps)) \\
&= (\ps\i)^*(\ze^M_X\o \ph\o\ps) = \ze^M_X\o\ph
     = \ze^P_X(\ph).
\endalign$$ 
The bundle projection 
$\ev_{x_0}:P_{x_0}(\g)\to M$ is visibly 
$\g$-equivariant.
Now we describe the associated unique principal connection form 
(Maurer-Cartan form) 
$\ka\in \Om^1(P_{x_0}(\g);\g)^{\Ga(\g)_{x_0}}$:
Consider a smooth curve $f_t$ in $P_{x_0}(\g)$. Then $\ddt f_t$ is a 
tangent vector with foot point $f_0$ and we have
$$
\ka(\ddt f_t) = (\operatorname{germ}_{x_0}\o \ze^M)\i 
     ((\ddt f_t)\o f_0\i)\in \g
$$

The dual action $\hat\ze^P:\g\to \X(P_{x_0}(\g))$ is given by 
$$
\hat\ze^P_X(\ph):= T\ph\o\ze^M_X = T\ph\o \ze^M_X\o\ph\i\o\ph = 
     ((\ph\i)^*\ze^M_X)\o\ph = \ze^M_{\Ad(\ph)X}\o\ph = 
     \ze^P_{\Ad(\ph)X}(\ph),
$$
by \nmb!{2.4} again, and its local flow is given by
$\Fl^{\hat\ze^P_X}_t(\ph) = 
\operatorname{germ}_{x_0}(\bar\ph\o\Fl^{\ze^M_X}_t)$, where $\bar\ph$ 
is a representative of the germ $\ph$. Then 
$\hat\ze^P:\g\to \X(P_{x_0}(\g))$. 
It is a Lie algebra anti homomorphism, since we have (using 
\cit!{9},~3.16)
$$\align
[\hat\ze^P_X,\hat\ze^P_Y](\ph) &= 
     \tfrac12\tfrac{\partial^2}{\partial t^2}|_0 
     \Fl^{\hat\ze^P_Y}_{-t} \Fl^{\hat\ze^P_X}_{-t} 
     \Fl^{\hat\ze^P_Y}_t \Fl^{\hat\ze^P_X}_t (\ph)\\
&= \tfrac12\tfrac{\partial^2}{\partial t^2}|_0 \ph \o \Fl^{\ze^M_X}_t 
     \o \Fl^{\ze^M_Y}_t\o \Fl^{\ze^M_X}_{-t}\o \Fl^{\ze^M_Y}_{-t}\\
&= T\ph\o [\ze^M_Y,\ze^M_X] = - T\ph\o \ze^M_{[X,Y]} = 
     -\hat\ze^P_{[X,Y]}(\ph).
\endalign$$
Consider now $\ps\in\Ga(\g)^{x_0}$ and its principal right action 
$\ps^*$ on $P_{x_0}(\g)$; it acts by pullback on each 
vector field $\hat\ze^P_X$ as follows:
$$\align
((\ps^*)^*\hat\ze^P_X)(\ph) &= (T(\ps^*)\i\o\hat\ze^P_X\o \ps^*)(\ph) 
     = (\ps\i)^*(\hat\ze^P_X(\ph\o\ps)) \\
&= (\ps\i)^*(T\ph\o T\ps\o \ze^M_X) = T\ph\o T\ps\o \ze^M_X\o \ps\i\\
&= T\ph\o ((\ps\i)^*\ze^M_X) = T\ph\o \ze^M_{\Ad(\ps)X} \\
&= \hat\ze^P_{\Ad(\ps)X}(\ph).
\endalign$$ 
Note that the vector fields $\hat\ze^P_{-X}\in\X(P_{x_0}(\g))$ for 
$X\in\g_{x_0}$ are the fundamental vector fields of the principal 
right action, and recall from \nmb!{2.6} that the Lie algebra of the 
structure group $\Ga(\g)_{x_0}$ is anti isomorphic to the isotropy Lie 
algebra $\g_{x_0}$.
The associated unique principal connection form 
(dual Maurer-Cartan form) 
$\hat\ka\in \Om^1(P_{x_0}(\g);\g)$ is given by 
$$\align
\hat\ka_\ph &= (\hat\ze^P_\ph)\i:T_\ph(P_{x_0}(\g)) \to \g \\
\hat\ka(\ddt f_t) &= (\operatorname{germ}_{x_0}\o\,  
     \ze^M)\i (Tf_0\i\o\ddt f_t)\in \g
\endalign$$
for each smooth curve $f_t$ in $P_{x_0}(\g)$. Since 
$\hat\ze^P:\g\to \X(P_{x_0}(\g))$ is a Lie algebra anti homomorphism, 
$\hat\ka$ satisfies the Maurer-Cartan equation in the form 
$d\hat\ka-\tfrac12[\hat\ka,\hat\ka]^\wedge=0$ and is 
$\Ga(\g)_{x_0}$-equivariant in the form $(\ps^*)^*\hat\ka = 
\Ad(\ps\i)\hat\ka$.

Finally let $\om\in\Om^1(M;\g)^\g$ be a principal connection form on 
the $\g$-manifold $M$. Then from \nmb!{5.7}.\thetag1 and from 
\nmb!{2.4} for any $\ph\in\Ga(\g)$ we have 
$$
\ph^*\om = \Ad(\ph\i)\om.
$$
We consider the 1-form 
$$\gather
\tilde\om_\ph := - \hat\ka_\ph - \Ad(\ph\i)\o\om_{\ph(x_0)}\o 
     T_\ph(\ev_{x_0}):T_\ph(P_{x_0}(\g))\to \g_{x_0}.\\
\tilde\om = - \hat\ka - (\Ad\o\operatorname{inv}).\ev_{x_0}^*\om.
\endgather$$
Then $\tilde\om$ is $\g_{x_0}$-valued by property 
\nmb!{5.7}.\thetag1, and it is a principal connection form
$\tilde\om\in\Om^1(P_{x_0}(\g);\g_{x_0})^{\Ga(\g)_{x_0}}$ on the principal 
$\Ga(\g)_{x_0}$-bundle $\ev_{x_0}:P_{x_0}(\g)\to M$, 
with right principal action now, $\ps^*(\ph)=\ph\o \ps$ for 
$\ps\in\Ga(\g)_{x_0}$ and $\ph\in P_{x_0}(\g)$, because we have in 
turn for $X\in\g_{x_0}$:
$$\align
\tilde\om_\ph(\ze^P_{-X}(\ph)) &= -\hat\ka_\ph(\ze^P_{-X}(\ph)) - 0 = X 
     \in \g_{x_0},\\
(\ps^*\tilde\om)_\ph &= \tilde\om_{\ph\o \ps}\o T(\ps^*) \\
&= - (\ps^*\hat\ka)_\ph - 
     \Ad(\ps\i\o\ph\i)\o\om_{\ph(x_0)}\o T(\ev_{x_0})\o T(\ps^*)\\
&= \Ad(\ps\i)(-\hat\ka_\ph - 
     \Ad(\ph\i)\o\om_{\ph(x_0)}\o T(\ev_{x_0}))\\
&= \Ad(\ps\i)\o\tilde\om_\ph.
\endalign$$
On the other hand, for any principal connection form 
$\tilde\om\in\Om^1(P_{x_0}(\g);\g_{x_0})^{\Ga(\g)_{x_0}}$ on the principal 
$\Ga(\g)_{x_0}$-bundle $\ev_{x_0}:P_{x_0}(\g)\to M$ the 
$\g$-valued 1-form $\bar\om_\ph:= \Ad(\ph)(-\hat\ka_\ph-\tilde\om_\ph)$ is 
horizontal and $\Ga(\g)_{x_0}$-invariant, thus it is the pullback of 
a unique form $\om\in\Om^1(M;\g)$ which is easily seen to be a 
principal connection form on $M$.
For the curvature we may compute as follows (compare \nmb!{5.6})
$$\align
\tilde\Om:&= d\tilde\om+\tfrac12[\tilde\om,\tilde\om]^\wedge_{\g_{x_0}} 
= d(- \hat\ka - (\Ad\o\operatorname{inv}).\ev_{x_0}^*\om) \\
&\quad  +\tfrac12[-\hat\ka-
     (\Ad\o\operatorname{inv}).\ev_{x_0}^*\om, - \hat\ka - 
     (\Ad\o\operatorname{inv}).\ev_{x_0}^*\om]
     ^\wedge_{\g_{x_0}} \\
&= -d\hat\ka+\tfrac12[\hat\ka,\hat\ka]^\wedge_\g 
     - d(\Ad\o\operatorname{inv})\wedge \ev_{x_0}^*\om - 
     (\Ad\o\operatorname{inv}).\ev_{x_0}^*d\om \\
&\quad - [\hat\ka,
     (\Ad\o\operatorname{inv}).\ev_{x_0}^*\om]^\wedge_{\g} 
     - \tfrac12[(\Ad\o\operatorname{inv}).\ev_{x_0}^*\om, 
     (\Ad\o\operatorname{inv}).\ev_{x_0}^*\om]^\wedge_\g\\
&=-(\Ad\o \operatorname{inv}).\ev_{x_0}^*(d\om + \tfrac12[\om,\om]^\wedge_\g) = 
     -(\Ad\o \operatorname{inv}).\ev_{x_0}^*\Om,
\endalign$$
where we used the Maurer-Cartan equation for $\hat\ka$ and 
$$\align
d(\Ad\o\operatorname{inv})(\hat\ze_X(\ph)) &= 
     \ddt (\Ad\o\operatorname{inv})(\ph\o\Fl^{\ze^M_X}_t) 
     = \ddt \Ad(\Fl^{\ze^M_X}_{-t}\o \ph\i) \\
&= -\ad(X)\o\Ad(\ph\i) = 
     -\ad(\hat\ka(\hat\ze^P_X(\ph)))\o\Ad(\ph\i)\\
d(\Ad\o\operatorname{inv}) &= 
-\ad(\hat\ka)(\Ad\o\operatorname{inv}).\qed
\endalign$$
\enddemo

\subhead\nmb.{5.9}. The Lie algebra $Z_{\X(M)}(\g)$ of infinitesimal 
automorphisms of a homogeneous $\g$-manifold \endsubhead
Let $M$ be an effective homogeneous $\g$-manifold. 
We will describe now the centralizer
$$
Z_{\X(M)}(\g):=\{\et\in\X(M):[\et,\ze_X]=0 \text{ for all }X\in\g\}
$$
of $\ze(\g)$ in the Lie algebra $\X(M)$ of all vector fields on $M$.

Let $x_0\in M$ be a fixed point with isotropy subalgebra 
$\g_{x_0}=\ker(\ze_{x_0}:\g\to T_{x_0}M)$ 
and isotropy group $\Ga(\g)_{x_0}$. 
We consider the \idx{\it normalizer} $N_\g(\Ga(\g)_{x_0})$ of the 
isotropy group $\Ga(\g)_{x_0}$ in $\g$, and the `Weyl algebra' $\hat\g = 
\hat\g(x_0)$, which are given by
$$\align
N_\g(\Ga(\g)_{x_0}) :&=\{X\in\g:\Ad(\ps)X-X\in\g_{x_0}
     \text{ for all }\ps\in\Ga(\g)_{x_0}\}\\
\hat\g(x_0) :&= N_\g(\Ga(\g)_{x_0})/\g_{x_0}.
\endalign$$
It is clear that $\g_{x_0}$ is an ideal in $N_\g(\Ga(\g)_{x_0})$, thus  
$\hat\g(x_0)$ is a Lie algebra. Clearly, if $\ph\in\Ga(\g)$, then 
$\Ad(\ph):\g\to\g$ induces an isomorphism 
$\hat\g(\ph(x_0))\to \hat\g(x_0)$.
We can define a Lie algebra anti homomorphism
$$
\hat\ze^M:\hat\g(x_0)=N_\g(\Ga(\g)_{x_0})/\g_{x_0}\to \X(M),
$$
as follows:
Any point $x\in M$ is of the form $x=\ph(x_0)$ for some element of the 
pseudogroup $\Ga(\g)$, and for
$X\in N_\g(\Ga(\g)_{x_0})$
we have a well defined vector field 
$$
\hat\ze^M_X(x=\ph(x_0)) := T\ph.\ze^M_X(x_0) = 
     ((\ph\i)^*\ze_X)(\ph(x_0)) = \ze_{\Ad(\ph)X}(x).
$$
For all $X\in \hat\g(x_0)$ and all $Y\in \g$ we have 
$$\align
((\Fl^{\ze_Y}_t)^*\hat\ze_X)(x=\ph(x_0)) :&= 
     T(\Fl^{\ze_Y}_{-t}).\hat\ze_X(\Fl^{\ze_Y}_t(\ph(x_0)))\\
&= T(\Fl^{\ze_Y}_{-t}).T(\Fl^{\ze_Y}_t\o\ph).\ze_X(x_0)\\
&= T(\ph).\ze_X(x_0) = \hat\ze_X(x=\ph(x_0)),\\
[\hat\ze_X,\ze_Y] &= 0,
\endalign$$
so that $\hat\ze(\hat\g(x_0))\subset \X(M)$ is contained in 
the centralizer $Z_{\X(M)}(\g)$. On the other hand we have:

\proclaim{Lemma} 
In this situation, $Z_{\X(M)}(\g) = \hat\ze(\hat\g)$, and 
these are exactly the vector fields on $M$ which are projections from 
all projectable vector fields in 
$\hat\ze^P(\g)\subset\X(P_{x_0}(\g))$ for the principal fiber bundle 
projection $\ev_{x_0}: P_{x_0}(\g)\to M$.
The flow of $\hat\ze_X$ for $X\in\hat\g(x_0)$ is given by 
$$
\Fl^{\hat\ze_X}_t(x=\ph(x_0)) = \ph(\Fl^{\ze_X}_t(x_0)).
$$
\endproclaim

\demo{Proof}
Let $\xi\in\X(M)$ be a vector field that commutes with the action of 
$\g$. Then for any $\ph\in\Ga(\g)$ we have 
$\ph^*\xi=T\ph\i\o\xi\o\ph=\xi$.
Then for $\ps\in \Ga(\g)_{x_0}$ we have $T_{x_0}\ps.\xi(x_0)= 
\xi_{x_0}$. If we choose any $X\in\g$ with $\ze_X(x_0)=\xi(x_0)$ we 
get by \nmb!{2.4} for all $\ps\in \Ga(\g)_{x_0}$
$$
\ze_{\Ad(\ps\i)X}(x_0)=(\ps^*\ze_X)(x_0) = T_{x_0}\ps.\xi(x_0) = 
\xi(x_0) = \ze_X(x_0).
$$
so that $X-\Ad(\ps\i)X\in\g_{x_0}$ and $X\in N_\g(\Ga(\g)_{x_0})$.
Moreover for $x\in M$ and $\ph\in\Ga(\g)$ with $x=\ph(x_0)$ we have 
$$
\hat\ze_X(x)= T_{x_0}\ph.\ze_X(x_0) = T_{x_0}\ph.\xi(x_0) = 
     \xi(\ph(x_0)) = \xi(x).
$$
The statement about the projectable vector fields on $P_{x_0(\g)}$ is 
easily checked, and the formula for the flow of $\hat\ze_X$ also 
follows by projecting it from $P_{x_0}(\g)$.
\qed\enddemo

\head\totoc\nmb0{6}. Parallel transport \endhead

\subhead\nmb.{6.1}. Local description of principal connections \endsubhead
Let $M$ be a locally trivial $\g$-manifold with projection 
$p:M\to N:=M/\g$ and with standard fiber $S$. 

Let $(U_\al,\ph_\al:p\i(U_\al)\to U_\al\x S)$ be an atlas of bundle 
charts as specified in \nmb!{2.1}.\therosteritem7. Then we have 
$(\ph_\be\o\ph_\al\i)(x,s)=(x,\ph_{\be\al}(x,s))$ for 
$(x,s)\in (U_\al\cap U_\be)\x S$, where $\ph_{\al\be}(x,\quad)$ is a 
$\g$-equivariant diffeomorphism of $S$ for each $x\in M$. See also 
\nmb!{5.9}.

Let $\Ph\in\Om^1(M;TM)^\g$ be a principal connection. Then we have
$$
((\ph_\al)^{-1})^*\Ph)(\xi_x,\et_y) =: 
     - \Ga^\al(\xi_x,y)  + \et_y \text{ for }\xi_x 
     \in T_xU_\al \text{ and } \et_y \in T_yS,
$$
since it reproduces vertical vectors. The $\Ga^\al$ are given by
$$
(0_x,\Ga^\al(\xi_x,y)) := - 
     T(\ph_\al).\Ph.T(\ph_\al)^{-1}.(\xi_x,0_y).
$$
We may consider $\Ga^\al$ as an element of the space 
$\Om^1(U_\al;Z_{\X(S)}(\g))$, i\.e\. as a
1-form on $U^\al$ with values in the centralizer $Z_{\X(S)}(\g)$
of $\ze^S(\g)$ in the Lie algebra $\X(S)$ of all vector fields on the 
standard fiber.  This space is finite dimensional by lemma \nmb!{5.9}. 
This follows from the naturality of the 
Fr\"olicher-Nijenhuis bracket \cit!{9},~8.15 via the following 
computation, with some abuse of notation:
$$
0=[\ze_X^M,\Ph]=\ph_\al^*[0\x \ze^S_X,\operatorname{pr}_2 - \Ga^\al]
     =\ph_\al^*(0\x[\ze^S_X,Id_{TS} - \Ga^\al])
     =-\ph_\al^*(0\x[\ze^S_X,\Ga^\al]),
$$
since $Id_{TS}$ is in the center of the Fr\"olicher-Nijenhuis 
algebra.

The $\Ga^\al$ are called the \idx{\it Christoffel forms} of the
connection $\Ph$ with respect to the bundle atlas $(U_\al,\ph_\al)$.

 From \cit!{9},~9.7 we get that the transformation law for the 
Christoffel forms is
$$
T_y(\ph_{\al\be}(x,\quad)).\Ga^\be(\xi_x,y) =
     \Ga^\al(\xi_x,\ph_{\al\be}(x,y)) 
     - T_x(\ph_{\al\be}(\quad,y)).\xi_x.
$$
The curvature $R$ of $\Ph$ satisfies 
$$
(\ph_\al^{-1})^*R = d\Ga^\al + \tfrac12[\Ga^\al,\Ga^\al]^\wedge_{\X(S)}.
$$
Here $d\Ga^\al$ is the exterior derivative of the 1-form 
$\Ga^\al\in\Om^1(U_\al;Z_{\X(S)}(\g))$ with values in the finite 
dimensional Lie algebra $Z_{\X(S)}(\g)$.

\subhead\nmb.{6.2}. Asystatic locally trivial $\g$-manifolds \endsubhead
A locally trivial $\g$-manifold $M$ is called asystatic if the 
normalizer $N_\g(\g_x)=\g_x$ for any isotropy subalgebra $\g_x$ of $M$. 
 From \nmb!{6.1} and \nmb!{5.9} we have immediately:

\proclaim{Proposition} An asystatic locally trivial $\g$-manifold admits a 
unique principal connection. This principal connection is flat. Its 
horizontal space at $x\in M$ is the subspace of $T_xM$ on which the 
isotropy representation of $\g_x$ vanishes.
\endproclaim 

\subhead\nmb.{6.3}. Horizontal lifts on locally trivial 
$\g$-manifolds \endsubhead
Let $\Ph$ be a connection on the locally trivial $\g$-manifold 
$(M,p,N=M/\g,S)$. Then the projection 
$(\pi_M,Tp):TM\to M\x_NTN$ onto the fibered product restricts to an 
isomorphism $\ker(\Ph)\to M\x_NTN$ whose inverse will be denoted by 
$C:M\x_NTN\to TM$ and will be called the \idx{\it horizontal lift}. 
If $\xi\in\X(N)$ is a vector field on the base then its horizontal 
lift $C(\xi)$ is given by $C(\xi)(y)=C(y,\xi(p(y)))$. In a bundle 
chart we have
$T(\ph_\al).C(\xi)(\ph_\al\i(x,s))= (\xi(x),\Ga^\al(\xi(x))(s))$.
Thus we see from \nmb!{6.1} that $C(\xi)$ commutes with all 
fundamental vector fields: $[\ze^M_X,C(\xi)]=0$ for all $X\in\g$ and 
$\xi\in\X(N)$.

Note that the $\g$-equivariant vector fields on $M$ which are 
horizontal in the sense that they take values in the kernel of the 
connection $\Ph$ are exactly the horizontal lifts of vector fields on 
the base manifold $N$.

\proclaim{\nmb.{6.4}. Theorem (Parallel transport)} Let $\Ph$ be a
connection on a locally trivial $\g$-manifold $(M,p,N=M/\g,S)$ 
and let $c:(a,b)\to N$ be
a smooth curve with $0 \in (a,b)$, $c(0) = x$. 

Then there is a neighborhood $U$ of $M_x \times \{0\}$ in
$M_x\times (a,b)$ and a smooth mapping $\Pt_c: U \to  M$ such
that:
\roster
\item $p(\Pt(c,u_x,t)) = c(t)$ if defined, and
$\Pt(c,u_x,0) = u_x$.
\item $\Ph(\frac d{dt}\Pt(c,u_x,t)) = 0$ if defined.
\item Reparametrization invariance: If $f: (a',b') \to   (a,b)$ is
smooth with  $0 \in (a',b')$, then 
$\Pt(c,u_x,f(t)) = \Pt(c \o f, \Pt(c,u_x,f(0)),t)$ 
if defined.    
\item $U$ is maximal for properties \therosteritem1 and \therosteritem2.
\item If the curve $c$ depends smoothly on further parameters then 
      $\Pt(c,u_x,t)$ depends also smoothly on those parameters.
\item If $\xi\in\X(N)$ is a vector field on the base and 
       $C(\xi)\in\X(M)$ is its horizontal lift, then 
       $\Pt(\Fl^\xi(x),u_x,t)=\Fl^{C(\xi)}_t(u_x)$.
\item For each $X\in \g$ the restrictions of the fundamental field 
       $\ze_X$ to $M_x=p\i(x)$ and to $M_{c(t)}$ are 
       $\Pt(c,t)$-related: 
       $T(\Pt(c,t))\o\ze_X|M_x=(\ze_X|M_{c(t)})\o \Pt(c,t)$.
\endroster
\endproclaim

\demo{Proof}
All assertions but the last two of this theorem follow from the 
general result \cit!{9},~9.8. The assertion \therosteritem6 is 
obvious and for \therosteritem7 we first note that it suffices to 
show it for curves of the form $c(t)=\Fl^\xi_t(x)$. But then by 
\therosteritem6 and by \nmb!{6.3} we have 
$$
\tfrac d{dt} \Pt(c,t)^*(\ze_X|M_{c(t)}) 
     =\tfrac d{dt} (\Fl^{C(\xi)}_t)^*(\ze_X)|M_x 
     =(\Fl^{C(\xi)}_t)^*([C(\xi),\ze_X])|M_x = 0
$$
so that $\Pt(c,t)^*(\ze_X|M_{c(t)})$ is constant in $t$ and thus 
equals $\ze_X|M_x$.
\qed\enddemo

\subhead\nmb.{6.5}. Parallel transport \endsubhead
Now we consider a $\g$-manifold $M$ which admits a principal connection 
$\Ph$. Guided by the last remark in \nmb!{6.3} we call 
\idx{\it parallel transport} each local flow $\Fl^\xi_t$ along any 
horizontal $\g$-equivariant vector field on $M$. 

\subhead\nmb.{6.6}. Complete connections \endsubhead
Let $M$ be a locally trivial $\g$-manifold with projection 
$p:M\to N:=M/\g$ and with standard fiber $S$. Following 
\cit!{9},~9.9 we call a principal connection $\Ph$ \idx{\it complete} 
if for each curve $c:(a,b)\to N$ the parallel transport $\Pt(c,\quad)$ is 
defined on the whole of $p\i(c(0))\x (a,b)$.

\proclaim{Proposition}
In this situation, if each vector field in the centralizer 
$Z_{\X(S)}(\g)$ of the $\g$-action on $S$ is complete, then each 
principal connection on $M$ is complete.
\endproclaim

\demo{Proof}
It suffices to show that for each curve $c:(a,b)\to U$ the parallel 
transport $\Pt(c,t)$ is defined on the whole of $M_{c(0)}$ for each 
$t\in(a,b)$, where $(U,\ph:M|U\to U\x S)$ is a bundle chart, since we 
may piece together such local solutions. So we 
may assume that $M=N\x S$ is a trivial $\g$-manifold. Then by 
\nmb!{6.1} any principal connection is of the form 
$\Ph(\xi_x,\et_s)=\et_s-\Ga(\xi_x)(s)$, where 
$\Ga\in\Om^1(N;\h)$ is the Christoffel form with values in 
the centralizer algebra $\h:=Z_{\X(S)}(\g)$, which is finite dimensional by 
\nmb!{5.9}. 
Since all vector fields in this Lie algebra are complete we may 
integrate its action on $S$ to a right action $r:S\x H\to S$ of a 
connected Lie group $H$ with Lie algebra 
$\h$. Then $t\mapsto \Ga(c'(t))$ is a smooth curve in 
$\h$ which we may integrate to a smooth curve $b(t)\in H$ 
with $b(0)=e$ and $b'(t)=L_{\Ga(c'(t))}(b(t))$ where $L_X$ is the 
left invariant vector field on $H$ generated by $X\in \h$.  It is an 
integral curve of a time dependent vector field on $H$ which is, 
locally in time, bounded with respect to a left invariant Riemannian 
metric on $H$. So indeed $b:(a,b)\to H$. 
But then $Pt(c,t,u)=(c(t), r(u,b(t)))$ for each $u\in M_{c(0)}$.
\qed\enddemo

\proclaim{\nmb.{6.7}.  Theorem}
Let $p:M\to N:=M/\g$ be a locally trivial $\g$-manifold with 
standard fiber $S$. Let $\Ph$ be a complete principal connection on 
$M$. Let us assume that the holonomy Lie algebra of $\Ph$ in the 
sense explained in the proof consists of complete vector fields on 
$S$.

Then there exists a finite dimensional Lie group $H$ with Lie 
algebra $\h$, a principal $H$-bundle $P\to N$, an irreducible 
principal connection form $\om$ on $P$, and a left action of 
$H$ on the standard fiber $S$ such that:
\roster
\item The fundamental vector field mapping of the $H$-action on $S$ 
       is an injective Lie algebra anti homomorphism 
       $\h\to Z_{\X(S)}(\g)$. The $H$-action on $S$ commutes 
       pointwise with the $\g$-action. 
\item The associated bundle $P[S]= P\x_HS$ is isomorphic to the 
       bundle $M\to N$.
\item The $\g$-principal connection $\Ph$ on $M$ is induced by the 
       principal connection form $\om$ on $P$.
\endroster 
\endproclaim

\demo{Proof} We suppose first that the base $N$ is connected. Let 
$x_0$ be a fixed point in $N$, and let us identify the standard fiber 
$S$ with the fiber $M_{x_0}$ of $M$ over $x_0$.
Since $\Ph$ is a complete connection on the bundle $M\to N$ we 
may consider the holonomy group $\operatorname{Hol}(\Ph,x_0)$ 
consisting of all parallel transports with respect to $\Ph$ along 
closed loops in $N$ through $x_0$, and the holonomy Lie algebra 
$\operatorname{hol}(\Ph,x_0)$, which is defined as follows (see 
\cit!{9},~9.10):

Let $C:TN\x_NM\to TM$ be the horizontal lift 
and let $R$ be the curvature of the connection $\Ph$.
For any $x\in N$ and $X_x\in T_xN$ the horizontal lift 
$C(X_x):= C(X_x,\quad): M_x\to TM$ is a vector field along
$M_x$. For $X_x$ and $Y_x\in T_xN$ we consider $R(CX_x,CY_x) \in
\X(M_x)$. Now we choose any piecewise smooth curve $c$ from
$x_0$ to $x$ and consider the diffeomorphism
$\Pt(c,t):S=M_{x_0}\to M_x$ and the pullback
$\Pt(c,1)^*R(CX_x,CY_x)\in \X(S)$. Then 
$\operatorname{hol}(\Ph,x_0)$ is the closed linear subspace,
generated by all these vector fields (for all $x\in N$, $X_x$,
$Y_x\in T_xN$ and curves $c$ from $x_0$ to $x$) in $\X(S)$
with respect to the compact $C^\infty$-topology.

In each local chart $(U_\al,\ph_\al:M|U\to U\x S)$ the curvature is 
expressed by the Christoffel form via
$(\ph_\al^{-1})^*R = d\Ga^\al + \tfrac12[\Ga^\al,\Ga^\al]^\wedge_{\X(S)}$, 
see \nmb!{6.1}, and since $\Ga^\al$ takes values in $Z_{\X(S)}(\g)$, the 
local expression of the curvature $(\ph_\al^{-1})^*R$  does it also.
The parallel transport $\Pt^\Ph(c,t)$ along any curves relates 
$\g$-fundamental vector fields to itself by \nmb!{6.4},~\thetag7.
Thus the holonomy Lie algebra $\operatorname{hol}(\Ph,x_0)$ is 
contained in the centralizer algebra $Z_{\X(S)}(\g)$, so it is finite 
dimensional.

By assumption $\operatorname{hol}(\Ph,x_0)\subset \X(S)$ consists of 
complete vector fields. Thus all conditions of theorem \cit!{9},~9.11 
are satisfied and all conclusions follow from it.
\qed\enddemo

\subhead\nmb.{6.8}. Remark \endsubhead
In the situation  of theorem \nmb!{6.7} let us suppose that the 
centralizer algebra $Z_{\X(S)}(\g)$ consists of complete vector 
fields. Then the each principal connection $\Ph$ is complete by 
\nmb!{6.6} and the holonomy Lie algebra 
$\operatorname{hol}(\Ph,x_0)\subset Z_{\X(S)}(\g)$ is also complete, 
see the proof of \nmb!{6.7}. Thus the conclusions of theorem 
\nmb!{6.7} hold.
                                                                 
\head\totoc\nmb0{7}. Characteristic classes for $\g$-manifolds \endhead

\subhead\nmb.{7.1}. Basic cohomology \endsubhead
Let $M$ be a $\g$-manifold. Following \nmb!{4.2},
by $\Om_{\text{hor}}^p(M)^\g$ we denote the space of all real valued 
horizontal forms on $M$ which are $\g$-invariant: $\L_{\ze_X}\ph=0$ 
for all $X\in\g$. These forms are called \idx{\it basic forms} of the 
$\g$-manifold $M$.

\proclaim{Lemma} In this situation the exterior derivative restricts 
to a mapping 
$$
d:\Om_{\text{hor}}^p(M)^\g\to\Om_{\text{hor}}^{p+1}(M)^\g
$$
\endproclaim

\demo{Proof}
Let $\ph\in\Om_{\text{hor}}^p(M)^\g$ then for $X\in\g$ we have 
$$\align
i_{\ze_X}d\ph &= i_{\ze_X}d\ph + di_{\ze_X}\ph =\L_{\ze_X}\ph =0\\
\L_{\ze_X}d\ph &= d\L_{\ze_X}\ph = 0. \qed
\endalign$$
\enddemo

The cohomology of the resulting differential complex will be called 
the \idx{\it basic cohomology} of the $\g$-manifold $M$:
$$
H^p_\g(M):= \frac
     {\ker(d:\Om_{\text{hor}}^p(M)^\g\to\Om_{\text{hor}}^{p+1}(M)^\g)}
     {\operatorname{im}(d:\Om_{\text{hor}}^{p-1}(M)^\g
     \to\Om_{\text{hor}}^p(M)^\g)}
$$
In the case of a $\g$-manifold $M$ of constant rank this cohomology 
is exactly the basic cohomology of the orbit foliation of $M$, 
defined by Reinhard \cit!{17} and intensively studied in the theory 
of foliations, see \cit!{13}, appendix~B by V\.~Sergiescu.
Note that this cohomology may be of infinite dimension, see \cit!{18} 
and \cit!{6}.

If $f:M\to N$ is a smooth $\g$-equivariant mapping between 
$\g$-manifolds $M$ and $N$, then the pullback operator induces a 
mapping $f^*:\Om_{\text{hor}}^p(N)^\g \to \Om_{\text{hor}}^p(M)^\g$ 
which in turn induces a linear mapping in basic cohomology
$f^*:H^p_\g(N) \to H^p_\g(M)$. If $f,g:M\to N$ are smoothly 
homotopic through $\g$-equivariant mappings then they induce the same 
mapping in basic cohomology.

\subhead\nmb.{7.2}. Chern-Weil forms \endsubhead
If $f\in L^k(\g) := (\bigotimes^k\g^*)$ is a $k$-linear function 
on $\g$ and if $\ps_i\in\Om^{p_i}(M;\g)$ 
we can construct the following differential forms (see \nmb!{4.1}):
$$\gather
\ps_1\otimes_\wedge \dots \otimes_\wedge \ps_k 
	\in \Om^{p_1+\dots+p_k}(M;\g\otimes \dots \otimes \g),\\
f^{\ps_1,\dots,\ps_k}:= f\o (\ps_1\otimes_\wedge \dots \otimes_\wedge \ps_k ) 
	\in \Om^{p_1+\dots+p_k}(M).
\endgather$$
The exterior derivative of the latter one is clearly given by
$$\multline
d( f\o(\ps_1\otimes_\wedge\dots\otimes_\wedge\ps_k)) = 
	f\o d(\ps_1\otimes_\wedge\dots\otimes_\wedge\ps_k) =\\
= f\o \left( \tsize\sum_{i=1}^k(-1)^{p_1+\dots+p_{i-1}} 
  	\ps_1\otimes_\wedge\dots\otimes_\wedge d\ps_i\otimes_\wedge
	\dots\otimes_\wedge\ps_k \right).
\endmultline$$
Note that the form $f^{\ps_1,\dots,\ps_k}$ is basic, i\.e\. 
$\g$-invariant and 
horizontal, if all $\ps_i\in \Om^{p_i}_{\text{hor}}(M;\g)^\g$ and 
$f$ is invariant under the adjoint action of $\g$ on 
$\g$ ($f\in L^k(\g)^\g$) in the following sense:

\definition{\nmb.{7.3}. Definition}
Let $\rh:\g\to\g\frak l(V)$ be a representation of $\g$. 
$f\in L^k(V)$ is called \idx{\it $\g$-invariant} if 
$\sum_{i=1}^k f(v_1,\dots,\rh(X)v_i,\dots,v_k)=0$ for each $X\in\g$.
If $f$ is $\g$-invariant then we have for $\ps_i\in\Om^{p_i}(M;V)$ and any 
$\ph\in\Om^p(M;\g)$, by applying alternation:
$$
f\o \left( \tsize\sum_{i=1}^k(-1)^{(p_1+\dots+p_{i-1})p} 
  	\ps_1\otimes_\wedge\dots\otimes_\wedge \rh^\wedge(\ph)\ps_i
	\otimes_\wedge	\dots\otimes_\wedge\ps_k \right) = 0.
$$
\enddefinition

\proclaim{\nmb.{7.4}. Lemma} Let $M$ be a $\g$-manifold with a 
principal connection form $\om$ and let $\Om$ be its curvature form.
Let $f\in L^k(\g)^\g$ be $\g$-invariant under the adjoint action then 
the differential form 
$f^\Om:=f^{\Om,\dots,\Om}\in \Om^{2k}(M)^\g$ is a closed 
$\g$-invariant form. 

If moreover $M$ is a free $\g$-manifold, then 
$\Om$ and consequently $f^\Om$ are horizontal, so 
$f^\Om\in \Om_{\text{hor}}^{2p}(M)^\g$ is a closed basic form.
\endproclaim

\demo{Proof}
We have in turn by \nmb!{7.2} and the Bianchi identity \nmb!{4.6}
$$\align
df^\Om &= d(f\o (\Om\otimes_\wedge \dots\otimes_\wedge \Om))\\
&= f\o \left(\tsize\sum_{i=1}^k 
     \Om\otimes_\wedge \dots\otimes_\wedge d\Om 
     \otimes_\wedge \dots\otimes_\wedge \Om\right)\\
&= -f\o \left(\tsize\sum_{i=1}^k 
     \Om\otimes_\wedge \dots\otimes_\wedge [\om,\Om]^\wedge 
     \otimes_\wedge \dots\otimes_\wedge \Om\right)  \\
&= -f\o \left(\tsize\sum_{i=1}^k 
     \Om\otimes_\wedge \dots\otimes_\wedge \ad^\wedge(\om)\Om 
     \otimes_\wedge \dots\otimes_\wedge \Om\right)
\endalign$$
which is 0 by \nmb!{7.3}. The second statement is obvious. 
\qed\enddemo

\proclaim{\nmb.{7.5}. Proposition} Let $\om_0$ and $\om_1$ be two 
principal connection forms on the $\g$-manifold $M$
with curvature forms $\Om_0, 
\Om_1\in \Om^2(M;\g)^\g$, and let $f\in L^k(\g)^\g$. 
Then the cohomology 
classes of the two closed forms $f^{\Om_0}$ and $f^{\Om_1}$ in 
$H^{2k}(M)$ coincide.

If $M$ is a free $\g$-manifold then the  curvature forms $\Om_0, 
\Om_1$ are horizontal and define the same basic cohomology classes in 
$H^{2k}_\g(M)=H^{2k}(\Om^*_{\text{hor}}(M)^{\g})$.
\endproclaim

Thus for $f\in L^k(\g)^\g$ the cohomology class 
$[f^\Om]\in H^{2p}(M)$ depends only on $f$ and the $\g$-action and 
we call it a \idx{\it characteristic class} for the $\g$-action. 

If $M\to M/G$ is a principal $G$-bundle, thus $M$ a free $\g$-manifold, 
we have just reconstructed the usual Chern-Weil characteristic classes.

If $M$ is a homogeneous $\g$-manifold (e\.g\. a homogeneous 
$G$-manifold $H\backslash G$), by theorem \nmb!{5.8} these 
characteristic classes in $H^{2m}(M)$ are usual characteristic 
classes of the principal $\Ga(\g)_{x_0}$-bundle $P_{x_0}(\g)\to M$, 
but possibly not all of them: only those arising from invariant 
polynomials on $\g_{x_0}$ which are restrictions of invariant 
polynomials on $\g$ appear.

\demo{Proof}
For each $t\in\Bbb R$ we have a principal connection form 
$\om_t:=(1-t)\om_0+t\om_1$, and also consider its curvature 
$\Om_t:=d\om_t + \frac12[\om_t,\om_t]^\wedge$. 
Since 
$\partial_t\om_t=\om_1-\om_0$ we get
$$\align
\partial_t \Om_t &= d\partial_t\om_t + [\om_t,\partial_t\om_t]^\wedge \\
&= d(\om_1-\om_0) + [\om_t,\om_1-\om_0]^\wedge = 
     d_{\om_t}(\om_1-\om_0).
\endalign$$
Note that $d_{\om_t}(\om_1-\om_0)$ makes sense since 
$\om_1-\om_0\in\Om^p_{\text{hor}}(M;\g)^\g$.
We will also need the Bianchi identity 
$d_{\om_t}\Om_t=d\Om_t+[\om_t,\Om_t]^\wedge =0$, see \nmb!{4.6}. 
Since $\Om_t$ is a 2-form we may assume that $f$ is symmetric. Then we 
have in turn:
$$\align
\partial_t f^{\Om_t} &= \partial_t f(\Om_t,\dots,\Om_t) = 
     p.f(\partial_t\Om_t,\Om_t,\dots,\Om_t)\\
&= p.f(d_{\om_t}(\om_1-\om_0),\Om_t,\dots,\Om_t)\\
&= p.f(d_{\om_t}(\om_1-\om_0),\Om_t,\dots,\Om_t)
     -p\sum_{i=2}^pf(\om_1-\om_0,\Om_t,\dots,d_{\om_t}\Om_t,\dots,\Om_t)\\
&= p.f(d(\om_1-\om_0),\Om_t,\dots,\Om_t)
     -p\sum_{i=2}^pf(\om_1-\om_0,\Om_t,\dots,d\Om_t,\dots,\Om_t)\\
&\qquad+ p.f([\om_t,\om_1-\om_0]^\wedge ,\Om_t,\dots,\Om_t)\\
&\qquad- p\sum_{i=2}^pf(\om_1-\om_0,\Om_t,\dots,
     [\om_t,\Om_t]^\wedge ,\dots,\Om_t)\\
&= p.df(\om_1-\om_0,\Om_t,\dots,\Om_t),
\endalign$$
where we again used \nmb!{7.3} in the form
$$
0=f([\om_t,\om_1-\om_0]^\wedge ,\Om_t,\dots,\Om_t)
     -\sum_{i+2}^pf(\om_1-\om_0,\Om_t,\dots,
     [\om_t,\Om_t]^\wedge ,\dots,\Om_t).
$$
Since 
$T(f,\Om_t):=f(\om_1-\om_0,\Om_t,\dots,\Om_t)\in 
     \Om^{2p-1}_{\text{hor}}(M)^\g$,
the following form is exact in $(\Om_{\text{hor}}^*(M)^\g,d)$:
$$
f^{\Om_1} - f^{\Om_0} = \int_0^1 \partial_t f^{\Om_t} dt 
= \int_0^1 pdT(f,\Om_t)dt = d\left(\int_0^1p.T(f,\Om_t)\,dt\right).
\qed$$
\enddemo

\enddocument